\documentclass[9pt,amsfonts]{article}
\usepackage{epsf,pgf,etex,graphicx,amssymb,eucal,mathrsfs}
\usepackage{latexsym,amsmath,epsfig,epic,eepic,xcolor}
\usepackage[all]{xy}
\usepackage{wasysym,ntheorem}
\usepackage{epsf,graphicx,pgf,etex,amscd,pstricks}

\usepackage[hypertexnames=false,backref=page,pdftex,
 	pdfpagemode=UseNone,
 	breaklinks=true,
 	extension=pdf,
 	colorlinks=true,
 	linkcolor=blue,
 	citecolor=red,
 	urlcolor=blue,
 ]{hyperref}

\newcommand{\Z}{\mathbb{Z}}
\newcommand{\R}{\mathbb{R}}

\newcommand{\cat}{{\rm CAT}(0)}

\newcommand{\hgotch}{\hat{\mathfrak{h}}}
\newcommand{\lgotch}{\hat{\mathfrak{l}}}
\newcommand{\kgotch}{\hat{\mathfrak{k}}}

\newtheorem{theorem}{{\bf Theorem}}

\newtheorem{prop}[theorem]{Proposition}
\newtheorem{lemma}[theorem]{Lemma}
\newtheorem{rem}{Remark}
\newtheorem{defi}{Definition}

\newtheorem{coro}{Corollary}

\theoremnumbering{Alph}
\newtheorem{main}{{\bf  Theorem}}
\newtheorem{cor}[main]{Corollary}

\title{Cubulable K\"ahler groups}
\author{Thomas Delzant and Pierre Py}
\date{September 2016}

\begin{document}

\maketitle

\begin{abstract} We prove that a K\"ahler group which is cubulable, i.e. which acts properly discontinuously and cocompactly on a $\cat$ cubical complex, has a finite index subgroup isomorphic to a direct product of surface groups, possibly with a free Abelian factor. Similarly, we prove that a closed aspherical K\"ahler manifold with a cubulable fundamental group has a finite cover which is biholomorphic to a topologically trivial principal torus bundle over a product of Riemann surfaces. Along the way, we prove a factorization result for essential actions of K\"ahler groups on irreducible, locally finite $\cat$ cubical complexes, under the assumption that there is no fixed point in the visual boundary.  
\end{abstract}

\tableofcontents

\newpage


\section{Introduction}

The purpose of this text is to characterize fundamental groups of compact K\"ahler manifolds which are {\it cubulable}, i.e. which act properly discontinuously and cocompactly on a $\cat$ {\it cubical complex}. For a survey concerning fundamental groups of compact K\"ahler manifolds, referred to as {\it K\"ahler groups} below, we refer the reader to~\cite{abckt,burger}. For basic definitions about cubical complexes, we refer the reader to~\cite[I.7]{bh} or~\cite{sageev1}. We only mention here that cubical complexes form a particular type of polyhedral complexes, and provide an important source of examples in geometric group theory.

From now on we will only deal with finite dimensional cubical complexes. Such complexes have a natural complete and geodesic metric. In the 80's, Gromov gave a combinatorial criterion for this metric to be $\cat$~\cite[\S 4.2]{gromov}. His work on {\it hyperbolization procedures}~\cite{gromov}, as well as Davis' work on reflection groups~\cite{davis} drew attention to these complexes; see also~\cite{basw,dj,djs}. Later, Sageev established a link between actions of a group $G$ on a $\cat$ cubical complex and the existence of a {\it multi-ended} subgroup of $G$~\cite{sageev0}. More recently, under the influence of Agol, Haglund and Wise, $\cat$ cubical complexes received a lot of attention as the list of examples of {\it cubulated groups} increased dramatically and thanks to their applications to $3$-dimensional topology. We refer the reader to~\cite{agol,bhw,bwise,hagenwise,haglundwise,wise0,wise1} for some of these developments. 

On the other hand, as we will see, actions of K\"ahler groups on $\cat$ cubical complexes are very constrained. Even more, K\"ahler groups are in a sense orthogonal to groups acting on $\cat$ cubical complexes. The first results in this direction go back to the work of the first author, together with Gromov~\cite{delzantgromov}. Note that most of the results of~\cite{delzantgromov} are not formulated in terms of actions on $\cat$ cubical complex but can be interpreted in terms of such actions thanks to the work of Sageev~\cite{sageev0}. In~\cite{delzantgromov}, the authors studied codimension one subgroups of hyperbolic K\"ahler groups. More generally they studied {\it filtered ends}, or {\it cuts} of hyperbolic K\"ahler groups, see section~\ref{filterede} for the definitions. Some of their results were later generalized by Napier and Ramachandran~\cite{nr2}. Stated informally, the results of~\cite{delzantgromov} show that under suitable hyperbolicity assumptions, the presence of sufficiently many subgroups of a K\"ahler group $\Gamma$ whose numbers of filtered ends are greater than $2$ implies that a finite index subgroup $\Gamma_{0}$ of $\Gamma$ is a {\it subdirect product} of a certain number of surface groups with a free Abelian group. All along this text, the expression {\it surface group} stands for the fundamental group of a closed oriented hyperbolic surface. Recall that a subdirect product of groups $G_{1}, \ldots ,G_{m}$ is a subgroup of $G_{1}\times \cdots \times G_{m}$ which surjects onto each factor. So, in the previous statement, we mean that $\Gamma_{0}$ is a subdirect product of $G_{1},\ldots , G_{m}$ where all the $G_{i}$'s are surface groups, except possibly for one of them which could be free Abelian. 

There is also a particular class of actions of K\"ahler groups on $\cat$ cubical complexes which are easy to describe. These are the ones given by homomorphisms into right-angled Artin groups. For homomorphisms from K\"ahler groups to right-angled Artin groups, it is easy to obtain a factorization result through a subdirect product of surface groups, possibly with a free Abelian factor. This relies on the facts that right-angled Artin groups embed into Coxeter groups and that Coxeter groups act properly on the product of finitely many trees. This was observed by the second author in~\cite{py}. 

These results left open the question of describing actions of K\"ahler groups on general $\cat$ cubical complexes. We answer this question here assuming that the cubical complexes are locally finite and that the actions have no fixed point in the visual boundary. We will briefly discuss in section~\ref{ques} the need for these two hypothesis. From this, one deduces easily a description of cubulable K\"ahler groups. 

The following statements involve various more or less standard notions about cubical complexes: essential actions, irreducible complexes, visual boundaries... We recall all these definitions in section~\ref{cch}.

\begin{main}\label{fac-pv} Let $M$ be a closed K\"ahler manifold whose fundamental group $\pi_{1}(M)$ acts on a finite dimensional locally finite irreducible $\cat$ cubical complex $X$. Assume that:
\begin{itemize} 
\item the action is essential, 
\item $\pi_{1}(M)$ has no fixed point in the visual boundary of $X$,
\item $\pi_{1}(M)$ does not preserve a Euclidean flat in $X$. 
\end{itemize}
Then there exists a finite cover $M_{1}$ of $M$ which fibers: there exists a holomorphic map with connected fibers $F : M_{1}\to \Sigma$ to a hyperbolic Riemann surface such that  the induced map $F_{\ast} : \pi_{1}(M) \to \pi_{1}(\Sigma)$ is surjective. Moreover the kernel of the homomorphism $F_{\ast}$ acts elliptically on $X$, i.e. its fixed point set in $X$ is nonempty.   
\end{main}

Hence, up to restricting the action to a convex subset of $X$, the action factors through the fundamental group of a hyperbolic Riemann surface. Indeed, as a convex subspace one can take the fixed point set of the subgroup $Ker(F_{\ast})$; this is a subcomplex of the first cubic subdivision of $X$.


From Theorem~\ref{fac-pv}, and using results due to Caprace and Sageev~\cite{cs} and to Bridson, Howie, Miller and Short~\cite{bhms2002}, we deduce a characterization of cubulable K\"ahler groups. In the next three statements, we implicitly assume that the cubical complexes are locally finite. 

\begin{main}\label{fac-complet} Suppose that a K\"ahler group $\Gamma$ acts properly discontinuously and cocompactly on a $\cat$ cubical complex $X$. Let $X=X_{1}\times \cdots \times X_{r}$ be the decomposition of $X$ into irreducible factors. Assume moreover that the action of $\Gamma$ on $X$ is essential and that each $X_{i}$ is unbounded. Then $\Gamma$ has a finite index subgroup $\Gamma_{\ast}$ which is isomorphic to a direct product 
\begin{equation}\label{peiso}
\Gamma_{\ast}\simeq H_{1}\times \cdots \times H_{k}\times G_{1} \times \cdots \times G_{m},
\end{equation}
where each $H_{j}$ is isomorphic to $\Z$, each $G_{j}$ is a surface group, and $k+m=r$ is the number of irreducible factors of $X$. Moreover the isomorphism~\eqref{peiso} can be chosen in such a way that for each $i$, $X_{i}$ contains a convex $\Gamma_{\ast}$-invariant subset $Y_{i}$ on which the action factors through one of the projections $\Gamma_{\ast}\to G_{j}$ or $\Gamma_{\ast}\to H_{j}$. 
\end{main}

If a group $G$ acts properly discontinuously and cocompactly on a $\cat$ cubical complex $X$, one can always find a $G$-invariant subcomplex $Y$ inside $X$, on which the action is essential  and all of whose irreducible factors are unbounded. This follows from Corollary 6.4 in~\cite{cs}. This implies the following corollary: 

\begin{cor}\label{cckcc} Suppose that a K\"ahler group $\Gamma$ acts properly discontinuously and cocompactly on a $\cat$ cubical complex $X$. Then $\Gamma$ has a finite index subgroup which is isomorphic to a direct product of surface groups with a free Abelian group. The conclusion of Theorem~\ref{fac-complet} holds after replacing $X$ by a suitable invariant subcomplex. 
\end{cor}

In what follows, we will say that a closed manifold $M$ is {\it cubulable} if it has the same homotopy type as the quotient of a finite dimensional $\cat$ cubical complex by a properly discontinuous, free, cocompact group action. Using an argument going back to Siu~\cite{siu}, and the fact that a cubulable manifold is aspherical, we will finally prove the following theorem. 

\begin{main}\label{manifold} If a K\"ahler manifold $M$ is cubulable, it admits a finite cover which is biholomorphic to a topologically trivial principal torus bundle over a product of Riemann surfaces. If $M$ is algebraic, it admits a finite cover which is biholomorphic to a direct product of an Abelian variety with a product of Riemann surfaces. 
\end{main}

We mention that the factorization result in Theorem~\ref{fac-pv} does {\it not} rely on the use of harmonic maps unlike most of the factorization results for actions of K\"ahler groups on symmetric spaces, trees or more general buildings. Indeed, although there is a general theory of harmonic maps with values into $\cat$ spaces~\cite{gs,ks}, the crucial step
$$harmonic \Rightarrow pluriharmonic$$
namely the use of a Bochner formula to prove that a harmonic map from a K\"ahler manifold to a certain nonpositively curved space is automatically pluriharmonic is not available when the target space is a $\cat$ cubical complex. Thus, our proof follows a different scheme. The idea is to produce fibrations of our K\"ahler manifold $M$ over a Riemann surface, in such a way that the kernel $N\lhd \pi_{1} (M)$ of the homomorphism induced by the fibration preserves a certain hyperplane $\hat{\mathfrak{h}}$ of the cubical complex. Since $N$ is normal it will have to preserve each hyperplane of the $\pi_{1}(M)$-orbit of $\hat{\mathfrak{h}}$. From this we will deduce that $N$ actually acts elliptically on the $\cat$ cubical complex, i.e. has a fixed point. 

From the K\"ahler point of view, the new idea of this paper appears in Section~\ref{pf}. It can be summarized as follows. Consider an infinite covering $Y$ of a closed K\"ahler manifold. Given a proper pluriharmonic function $u$ on $Y$, one can look for conditions under which the foliation induced by the holomorphic form $du^{1,0}$ is actually a fibration. Our observation is that it is enough to find on $Y$ a second {\it plurisubharmonic} function which is not a function of $u$. This differs from other fibration criterions in the study of K\"ahler groups as the second function we need is only required to be plurisubharmonic. In the various Castelnuovo-de Franchis type criterions known so far, one needs two pluriharmonic functions, or two holomorphic $1$-forms, see~\cite{nrjta} and the references mentioned in the introduction there. 

The text is organized as follows. In section~\ref{gaoccc}, we recall basic facts on $\cat$ cubical complexes as well as more advanced results due to Caprace and Sageev. Given a group $G$ acting on a $\cat$ cubical complex and satisfying suitable hypothesis, wee also explain how to construct a hyperplane $\hgotch$ whose stabilizer $H$ in $G$ has the property that the Schreier graph of $H\backslash G$ is non-amenable. In section~\ref{rems}, we explain how to construct some nontrivial plurisubharmonic functions on certain covering spaces of a K\"ahler manifold whose fundamental group acts on a $\cat$ cubical complex. The proof of Theorem~\ref{fac-pv} is concluded in Section~\ref{pf}. In section~\ref{ckmag}, we prove Theorem~\ref{fac-complet} and~\ref{manifold}. Finally, section~\ref{ques} contains a few comments about possible improvements of our results.

{\bf Acknowledgements.} We would like to thank Pierre-Emmanuel Caprace, Yves de Cornulier and Misha Sageev for many explanations about $\cat$ cubical complexes as well as Martin Bridson for discussions motivating this work.

\section{Groups acting on ${\rm CAT}(0)$ cubical complexes}\label{gaoccc}

 In this section we recall some basic properties of $\cat$ cubical complexes and some advanced results due to Caprace and Sageev~\cite{cs}. For a more detailed introduction to these spaces, we refer the reader to~\cite[I.7]{bh} and \cite{sageev1} for instance. From now on and until the end of the text, all $\cat$ cubical complexes will be finite dimensional; we will not mention this hypothesis anymore. By convention, we assume that all cubes of our cubical complexes are isometric to $[-1,1]^{n}$ for some $n$. 

\subsection{Cubical complexes and hyperplanes}\label{cch}

A cubical complex $X$ can be naturally endowed with a distance as follows. If $x$ and $y$ are in $X$, one considers chains
$$x_{0}=x, x_{1}, \ldots , x_{n}=y$$
where each pair $(x_{i}, x_{i+1})$ is contained in a cube $C_{i}$ of $X$ and one defines 
$$d(x,y)={\rm inf} \sum_{i=0}^{n-1}d_{C_{i}}(x_{i},x_{i+1})$$
where the inf is taken over all possible chains. Here $d_{C_{i}}$ is the intrinsic distance of the cube; the number $d_{C_{i}}(x_{i},x_{i+1})$ does not depend on the choice of the cube containing $x_{i}$ and $x_{i+1}$. The function $d$ is a distance which is complete and geodesic thanks to the finite dimensionality hypothesis~\cite[I.7]{bh}. From now on, we will always assume that the cube complexes we consider are $\cat$ spaces, when endowed with the previous metric. According to a classical theorem by Gromov~\cite[\S 4.2]{gromov}, this is equivalent to the fact that the complex is simply connected and that the link of each vertex is a flag complex. The visual boundary~\cite[II.8]{bh} of a $\cat$ cubical complex $X$ is denoted by $\partial_{\infty}X$.

We now recall the definition of {\it hyperplanes}, first introduced by Sageev~\cite{sageev0}. We fix a $\cat$ cubical complex $X$. Let $\Box$ be the equivalence relation on edges of our complex,
defined as follows. One says that two edges $e$ and $f$ are equivalent, denoted $e\Box f$, if there exists a chain $$e_{1}=e, \ldots , e_{n}=f$$
such that for each $i$, $e_{i}$ and $e_{i+1}$ are opposite edges of some $2$-dimensional cube. If $e$ is an edge we will denote by $[e]$ its equivalence class. 

A {\it midcube} of a cube $C$, identified with $[-1,1]^{n}$, is the subset of $C$ defined by setting one of the coordinates equal to $0$. One says that a midcube and an edge $e$ of a cube $C \simeq [-1,1]^{n}$ are {\it transverse} if the midcube is
defined by $t_{i}=0$ and the coordinate $t_{i}$ is the only nonconstant coordinate of the egde $e$. Now the hyperplane associated to an equivalence class of edges $[e]$ is the union of all midcubes which are transverse to an edge of the class $[e]$. It is denoted by $\hat{\mathfrak{h}}([e])$. If we want to denote a hyperplane without referring to the edge used to define it, we will denote it by $\hat{\mathfrak{h}}$. Finally we will denote by $N(\hat{\mathfrak{h}})$ the union of the interiors of all cubes which intersect a hyperplane $\hat{\mathfrak{h}}$. 

One can prove that hyperplanes enjoy the following properties~\cite{sageev0}:
\begin{enumerate}
\item Each hyperplane intersects a cube in at most one midcube.
\item Each hyperplane $\hat{\mathfrak{h}}$ separates $X$ into two connected components, the closures of the two connected components of $X-\hat{\mathfrak{h}}$ are called the two {\it half-spaces} associated to $\hat{\mathfrak{h}}$.
\item Every hyperplane as well as every halfspace is a convex subset of $X$ for the distance $d$.
\item For every hyperplane $\hat{\mathfrak{h}}$, the set $N(\hat{\mathfrak{h}})$ is convex and naturally isometric to $\hat{\mathfrak{h}} \times (-1,1)$. If $t\in (-1,1)$, the subset of $N(\hgotch)$ corresponding to $\hgotch \times \{t\}$ under this isometry is called a {\it translated hyperplane} of $\hgotch$. 
\item If a geodesic segment intersects a hyperplane in two points, it is entirely contained in it. 
\end{enumerate}

The group ${\rm Aut}(X)$ of automorphisms of a $\cat$ cubical complex $X$ is the group of all permutations $X \to X$ which send isometrically a cube of $X$ to another cube of $X$. An automorphism of $X$ is automatically an isometry of the distance $d$ introduced earlier. In what follows, every time we will say that a group $G$ acts on a $\cat$ cubical complex, we will mean that we have a homomorphism from the group $G$ to the group ${\rm Aut}(X)$ of automorphisms of $X$. In this case, if $\hat{\mathfrak{h}}$ is a hyperplane, we will denote by 
$$Stab_{G}(\hat{\mathfrak{h}})$$
the subgroup of $G$ of elements which globally preserve $\hat{\mathfrak{h}}$ and which also preserve each component of the complement of $\hat{\mathfrak{h}}$. 

Following~\cite{cs}, we say that an action of a group $G$ on a $\cat$ cubical complex $X$ is {\it essential} if no $G$-orbit is contained in a bounded neighborhood of a halfspace. This implies in particular that $G$ has no fixed point in $X$. We will also use the following convention. If $\mathfrak{h}$ is a halfspace of $X$, we will denote by $\hat{\mathfrak{h}}$ the associated hyperplane (the boundary of $\mathfrak{h}$) and by $\mathfrak{h}^{\ast}$ the opposite halfspace, i.e. the closure of $X\setminus\mathfrak{h}$. 

Finally, we mention that there is a natural notion of irreducibility for a $\cat$ cubical complex, see~\cite[\S 2.5]{cs}. Any finite dimensional $\cat$ cube complex $X$ has a canonical decomposition as a product of finitely many irreducible $\cat$ cubical complexes. Moreover every automorphism of $X$ preserves this decomposition, up to a permutation of possibly isomorphic factors. We refer the reader to~\cite[\S 2.5]{cs} for a proof of these facts. We will use this decomposition to deduce Theorem~\ref{fac-complet} from Theorem~\ref{fac-pv}.

\subsection{Existence of strongly separated hyperplanes}\label{essh}

The following definition is due to Behrstock and Charney~\cite{berch} and is a key notion to study {\it rank $1$ phenomena} for group actions on $\cat$ cubical complexes.

\begin{defi} Two hyperplanes $\hat{\mathfrak{h}}_{1}$ and $\hat{\mathfrak{h}}_{2}$ in a $\cat$ cubical complex are strongly separated if they are disjoint and if there is no hyperplane meeting both $\hat{\mathfrak{h}}_{1}$ and $\hat{\mathfrak{h}}_{2}$.
\end{defi}

\noindent If $X$ is a $\cat$ cubical complex and $Y$ is a closed convex subset of $X$, we will denote by $\pi_{Y} : X \to Y$ the projection onto $Y$~\cite[II.2]{bh}. The following proposition is also taken from~\cite{berch}.

\begin{prop}\label{projunique}
If the hyperplanes $\hat{\mathfrak{h}}_{1}$ and $\hat{\mathfrak{h}}_{2}$ are strongly separated, there exists a unique pair of point $(p_{1}, p_{2}) \in \hat{\mathfrak{h}}_{1} \times \hat{\mathfrak{h}}_{2}$ such that $d(p_{1},p_{2})=d(\hat{\mathfrak{h}}_{1},\hat{\mathfrak{h}}_{2})$. The projection of any point of $\hat{\mathfrak{h}}_{2}$ (resp. $\hat{\mathfrak{h}}_{1}$) onto $\hat{\mathfrak{h}}_{1}$ (resp. $\hat{\mathfrak{h}}_{2}$) is $p_{1}$ (resp. $p_{2}$).
\end{prop}
{\it Proof.} The first claim is Lemma 2.2 in~\cite{berch}. The proof given there also shows that if $\hgotch$ is a hyperplane distinct from $\hat{\mathfrak{h}}_{1}$ and $\hat{\mathfrak{h}}_{2}$, no translated hyperplane of $\hgotch$ can intersect both $\hat{\mathfrak{h}}_{1}$ and $\hat{\mathfrak{h}}_{2}$. We now prove the last assertion of the proposition. It is enough to prove that the projection on $\hat{\mathfrak{h}}_{1}$ of each point of $\hat{\mathfrak{h}}_{2}$ is the middle of an edge of the cubical complex. Since $\hat{\mathfrak{h}}_{1}$ is connected, this implies that all points of $\hat{\mathfrak{h}}_{1}$ have the same projection, which must necessarily be the point $p_{1}$. Let $x\in \hat{\mathfrak{h}}_{2}$. If the projection $q$ of $x$ onto $\hat{\mathfrak{h}}_{1}$ is not the middle of an edge, than there exists a cube $C$ of dimension at least $3$ which contains $q$ as well as the germ at $q$ of the geodesic from $q$ to $x$. One can identify $C$ with $[-1,1]^{n}$ in such a way that 
$$q=(0,s,s_{3}, \ldots , s_{n})$$
with $\vert s \vert <1$. Since the germ of geodesic going from $q$ to $x$ is orthogonal to $\hat{\mathfrak{h}}_{1}$, it must be contained in $[-1,1]\times \{s\}\times [-1,1]^{n-2}$. We call $\hat{\mathfrak{m}}$ the hyperplane associated to any edge of $C$ parallel to $\{0\}\times [-1,1]\times \{0\}^{n-2}$. Hence the germ of $[q,x]$ at $q$ is contained in a translated hyperplane of $\hat{\mathfrak{m}}$. This implies that $[q,x]$ is entirely contained in this translated hyperplane. This contradicts the fact that no translated hyperplane can intersect both $\hat{\mathfrak{h}_{1}}$ and $\hat{\mathfrak{h}_{2}}$.\hfill $\Box$

\noindent Note that the second point of the proposition can be stated in the following slightly stronger way. If $\mathfrak{h}_{1}$ is the half-space defined by $\hat{\mathfrak{h}}_{1}$ and which does not contain $\hat{\mathfrak{h}}_{2}$, then the projection of any point of $\mathfrak{h}_{1}$  onto $\hat{\mathfrak{h}}_{2}$ is equal to $p_{2}$. Indeed if $q\in \mathfrak{h}_{1}$ and if $\gamma$ is the geodesic from $q$ to $\pi_{\hat{\mathfrak{h}}_{2}}(q)$, there must exist a point $q'$ on $\gamma$ which belongs to $\hat{\mathfrak{h}}_{1}$. Hence $\pi_{\hat{\mathfrak{h}}_{2}}(q)=\pi_{\hat{\mathfrak{h}}_{2}}(q')=p_{2}$. This proves the claim. 

The following theorem is due to Caprace and Sageev~\cite[\S 1.2]{cs}. It gives sufficient conditions for the existence of strongly separated hyperplanes in a $\cat$ cubical complex $X$. 

\begin{theorem}\label{existence-strong-s} Assume that the $\cat$ cubical complex is irreducible and that the group ${\rm Aut}(X)$ acts essentially on $X$ and without fixed point in the visual boundary. Then $X$ contains two strongly separated hyperplanes. 
\end{theorem}

In the end of this section, we consider a $\cat$ cubical complex $X$, two strongly separated hyperplanes $\hat{\mathfrak{h}}$ and $\hat{\mathfrak{k}}$ in $X$, and a group $G$ acting on $X$. We prove a few lemmas which will be used in the next sections. We denote by $\mathfrak{k}$ the halfspace delimited by $\hat{\mathfrak{k}}$ which does not contain $\hat{\mathfrak{h}}$. 

\begin{lemma} Let $p$ be the projection of $\hat{\mathfrak{k}}$ (or $\mathfrak{k}$) on $\hat{\mathfrak{h}}$. Let $h$ be an element of $Stab_{G}(\hat{\mathfrak{h}})$. If $h(\mathfrak{k})\cap \mathfrak{k}$ is nonempty, then $h$ fixes $p$. 
\end{lemma}
{\it Proof.} This is an easy consequence of Proposition~\ref{projunique}. Let $\pi_{\hat{\mathfrak{h}}} : X \to \hat{\mathfrak{h}}$ be the projection. We have seen that $\pi_{\hat{\mathfrak{h}}} (\mathfrak{k})=p$. Since $h\in Stab_{G}(\hat{\mathfrak{h}})$, the map $\pi_{\hat{\mathfrak{h}}}$ is $h$-equivariant. Let $x\in \mathfrak{k}$ be such that $h(x)\in \mathfrak{k}$. We have:
$$\pi_{\hat{\mathfrak{h}}}(h(x))=p$$
since $h(x)\in \mathfrak{k}$. But we also have $\pi_{\hat{\mathfrak{h}}}(h(x))=h(\pi_{\hat{\mathfrak{h}}}(x))=h(p)$ since $x\in \mathfrak{h}$. Hence $h(p)=p$.\hfill $\Box$

\medskip

We now define
$$\Sigma=\{h\in Stab_{G}(\hat{\mathfrak{h}}), h(\mathfrak{k})\cap \mathfrak{k} \neq \emptyset \}$$
and let $A$ be the subgroup of $Stab_{G}(\hat{\mathfrak{h}})$ generated by $\Sigma$. According to the previous lemma, every element of $A$ fixes $p$. Let
 $$U=\bigcup_{a\in A}a(\mathfrak{k})$$
be the union of all the translates by $A$ of the halfpsace $\mathfrak{k}$. 

\begin{lemma}\label{separea} Let $h$ be an element of $Stab_{G}(\hgotch)$. If $h(U)\cap U$ is nonempty, then $h\in A$ and $h(U)=U$. 
\end{lemma}
{\it Proof.} Let $h$ be such that $h(U)\cap U$ is nonempty. Then there exist $a_{1}$ and $a_{2}$ in $A$ such that $ha_{1}(\mathfrak{k})\cap a_{2}(\mathfrak{k})\neq \emptyset$. This implies that $a_{2}^{-1}ha_{1}$ is in $\Sigma$, hence in $A$. In particular $h$ is in $A$.\hfill $\Box$

If the space $X$ is proper, there are only finitely many hyperplanes passing through a given ball of $X$. Since the group $A$ fixes the point $p$, in this case we get that the family of hyperplanes 
$$\left( a(\kgotch )\right)_{a\in A}$$
is actually a {\it finite} family. This implies that the family of halfspaces $(a(\mathfrak{k}))_{a\in A}$ is also finite. Note that the same conclusion holds if we assume that the action of $G$ has finite stabilizers instead of assuming the properness of $X$. Indeed, the whole group $A$ is finite in that case. We record this observation in the following:

\begin{prop}\label{finitude} If $X$ is locally finite or if the $G$-action on $X$ has finite stabilizers, then there exists a finite set $A_{1}\subset A$ such that for all $a\in A$, there exists $a_{1}\in A_{1}$ such that $a_{1}(\mathfrak{k})=a(\mathfrak{k})$. 
\end{prop}

We will make an important use of this Proposition in section~\ref{sphm}. Although we will only use it under the hypothesis that $X$ is locally finite, we decided to remember the fact that the conclusion still holds for actions with finite stabilizers as this might be useful to study proper actions of K\"ahler groups on non-locally compact $\cat$ cubical complexes.

\subsection{Non-amenability of certain Schreier graphs}

In this section we consider an irreducible $\cat$ cubical complex and a finitely generated group $G$ which acts essentially on $X$, does not preserve any Euclidean subspace of $X$, and has no fixed point in $\partial_{\infty}X$. A consequence of the work of Caprace and Sageev~\cite{cs} is that under these hypothesis, the group $G$ contains a nonabelian free group. We will need the following slight modification of this important fact; see also~\cite{kasa} for a similar statement.

To state the next theorem, we need the following definition. Following Caprace and Sageev~\cite{cs}, we say that three halfspaces $\mathfrak{a}$, $\mathfrak{b}$, $\mathfrak{c}$ form a {\it facing triple of halfspaces} if they are pariwise disjoint.  

\begin{theorem}\label{groupeslibres} Let $X$ be an irreducible $\cat$ cubical complex. Assume that $G$  is a finitely generated group of automorphisms of $X$ which satisfies the following three conditions: $G$ does not fix any point in the visual boundary of $X$, does not preserve any Euclidean subspace of $X$ and acts essentially. 

Then $X$ contains a facing triple of halfspaces $\mathfrak{k}$, $\mathfrak{h}$, $\mathfrak{l}$ such that the three hyperplanes $\kgotch$, $\hgotch$, $\lgotch$ are strongly separated and such that there exists a non-Abelian free group $F<G$ with the property that $F\cap gStab_{G}(\hgotch)g^{-1} =\{1\}$ for all $g\in G$.  
\end{theorem}

Besides the facts about $\cat$ cubical complexes already recalled in the previous sections, we will further use the following three results from~\cite{cs}, which apply under the hypothesis of the previous theorem.
\begin{enumerate}
\item The space $X$ contains a facing triple of halfspaces, see Theorem E in~\cite{cs}.\label{exft} 
\item We will use the {\it flipping lemma} from~\cite[\S 1.2]{cs}: for any halfspace $\mathfrak{h}$, there exists $g\in G$ such that $\mathfrak{h}^{\ast} \subsetneq g(\mathfrak{h})$.
\item Finally we will also use the {\it double skewer lemma} from~\cite[\S 1.2]{cs}: for any two halfspaces $\mathfrak{k}\subset \mathfrak{h}$, there exists $g\in G$ such that $g(\mathfrak{h})\subsetneq \mathfrak{k} \subset \mathfrak{h}$.
\end{enumerate}

We now turn to the proof of the theorem.

\noindent {\it Proof of Theorem~\ref{groupeslibres}.} By~\ref{exft} above, one can choose a facing triple of halfspaces
$$\mathfrak{h}, \mathfrak{h}_{1}, \mathfrak{h}_{2}$$
in $X$. By the flipping lemma, there exists an element $k\in G$ such that $k (\mathfrak{h}^{\ast})\subsetneq \mathfrak{h}$. We now define $\mathfrak{h}_{3}=k(\mathfrak{h}_{1})$ and $\mathfrak{h}_{4}=k(\mathfrak{h}_{2})$. By construction, $\mathfrak{h}_{1}$, $\mathfrak{h}_{2}$, $\mathfrak{h}_{3}$, $\mathfrak{h}_{4}$ is a facing quadruple of halfspaces. We will need to assume moreover that these four halfspaces are strongly separated. This can be done thanks to the following lemma.

\begin{lemma}\label{interm} There exists half-spaces $\mathfrak{h}_{j}'\subset \mathfrak{h}_{j}$ ($1\le j \le 4$) such that the $\mathfrak{h}_{j}'$s are strongly separated. 
\end{lemma}
{\it Proof of Lemma~\ref{interm}.} According to Theorem~\ref{existence-strong-s}, we can find two halfspaces $\mathfrak{a}_{1}\subset \mathfrak{a}_{2}$ such that the corresponding hyperplanes $\hat{\mathfrak{a}}_{i}$ are strongly separated. We claim the following:
\begin{center}
{\it Up to replacing the pair $(\mathfrak{a}_{1},\mathfrak{a}_{2})$ by the pair $(\mathfrak{a}_{2}^{\ast},\mathfrak{a}_{1}^{\ast})$, there exists $i\in \{ 1, 2, 3, 4\}$  and $x\in G$ such that $x(\mathfrak{a}_{1})\subset x(\mathfrak{a}_{2})\subset \mathfrak{h}_{i}^{\ast}$.}
\end{center}
Let us prove this claim. First we prove that one of the four halfspaces $\mathfrak{a}_{1}$, $\mathfrak{a}_{2}$, $\mathfrak{a}_{1}^{\ast}$, $\mathfrak{a}_{2}^{\ast}$ is contained in $\mathfrak{h}_{j}^{\ast}$ for some $j$. If this is false, each of these four halfspaces intersects the interior of each of the $\mathfrak{h}_{j}$. In particular $\mathfrak{a}_{k}$ and $\mathfrak{a}_{k}^{\ast}$ intersect the interior of $\mathfrak{h}_{j}$. Since $\mathfrak{h}_{j}$ is convex, $\hat{\mathfrak{a}}_{k}$ intersects it also. Considering two indices $j\neq j'$ and a geodesic from a point in $\hat{\mathfrak{a}}_{k}\cap \mathfrak{h}_{j}$ to a point in $\hat{\mathfrak{a}}_{k}\cap \mathfrak{h}_{j'}$, one sees that $\hat{\mathfrak{a}}_{k}$ intersects each of the hyperplanes $\hat{\mathfrak{h}}_{j}$. Since this is true for $k=1, 2$, this contradicts the fact that the hyperplanes $\hat{\mathfrak{a}_{1}}$ and $\hat{\mathfrak{a}_{2}}$ are strongly separated. This concludes the proof that one of the four halfspaces $\mathfrak{a}_{1}$, $\mathfrak{a}_{2}$, $\mathfrak{a}_{1}^{\ast}$, $\mathfrak{a}_{2}^{\ast}$ is contained in $\mathfrak{h}_{j}^{\ast}$ for some $j$. If $\mathfrak{a}_{2}$ or $\mathfrak{a}_{1}^{\ast}$ is contained in one of the $\mathfrak{h}_{j}^{\ast}$, this proves the claim (with $x=1$). Otherwise we assume that $\mathfrak{a}_{1}\subset \mathfrak{h}_{j}^{\ast}$ (the last case $\mathfrak{a}_{2}^{\ast}\subset \mathfrak{h}_{j}^{\ast}$ being similar). In this case the double skewer lemma applied to $\mathfrak{a}_{1}\subset \mathfrak{a}_{2}$ implies that there exists $x\in G$ such that:
$$x(\mathfrak{a}_{1})\subset x(\mathfrak{a}_{2})\subset \mathfrak{a}_{1}\subset \mathfrak{h}_{j}^{\ast}.$$
This proves our claim. We now write $\mathfrak{b}_{1}:=x(\mathfrak{a}_{1})$, $\mathfrak{b}_{2}:=x(\mathfrak{a}_{2})$.

Since $\mathfrak{h}_{1}\subset \mathfrak{h}_{2}^{\ast}$, the double skewer lemma implies that there exists $g\in G$ such that $g(\mathfrak{h}_{2}^{\ast})\subsetneq \mathfrak{h}_{1}$. Similarly there exists $h\in G$ such that $\mathfrak{h}_{3} \supsetneq h(\mathfrak{h}_{4}^{\ast})$. Applying one of the four elements $g, g^{-1}, h, h^{-1}$ to $\mathfrak{b}_{1}$ and $\mathfrak{b}_{2}$, we obtain two halfspaces which are contained in one of the $\mathfrak{h}_{j}$. For instance if $i=2$ in the claim above, one has $g(\mathfrak{b}_{1})\subset g(\mathfrak{b}_{2})\subset \mathfrak{h}_{1}$. In what follows we continue to assume that we are in this case, the other ones being similar.  Since $\mathfrak{h}_{1}\subset \mathfrak{h}_{4}^{\ast}$ and since $h(\mathfrak{h}_{4}^{\ast})\subset \mathfrak{h}_{3}$, one has  
$$hg(\mathfrak{b}_{1})\subset hg(\mathfrak{b}_{2})\subset \mathfrak{h}_{3}.$$
Finally by a similar argument we have: 
$$g^{-1}hg(\mathfrak{b}_{1})\subset g^{-1}hg(\mathfrak{b}_{2})\subset \mathfrak{h}_{2}\;\;\;\; and \;\;\;\; h^{-1}g(\mathfrak{b}_{1})\subset h^{-1}g(\mathfrak{b}_{2})\subset \mathfrak{h}_{4}.$$ 
We now define $\mathfrak{h}_{1}'=g(\mathfrak{b}_{1})$, $\mathfrak{h}_{2}'=g^{-1}hg(\mathfrak{b}_{1})$, $\mathfrak{h}_{3}'=hg(\mathfrak{b}_{1})$ and $\mathfrak{h}_{4}'=h^{-1}g(\mathfrak{b}_{1})$. We check that the $\mathfrak{h}_{j}'$s are strongly separated halfspaces. It is clear that they are pairwise disjoint. We check that the corresponding hyperplanes are strongly separated. We do this for the pair  $\{ \hat{\mathfrak{h}_{1}'} , \hat{\mathfrak{h}_{2}'} \}$, the other cases being similar. If a hyperplane $\hat{\mathfrak{u}}$ intersects both $\hat{\mathfrak{h}_{1}'}$ and $\hat{\mathfrak{h}_{2}'}$, it will have to intersect also $g(\hat{\mathfrak{b}}_{2})$. This contradicts the fact that the pair 
$$g(\hat{\mathfrak{b}}_{1}), g(\hat{\mathfrak{b}}_{2})$$ is strongly separated. Hence $\hat{\mathfrak{h}_{1}'}$ and $\hat{\mathfrak{h}_{2}'}$ are strongly separated.\hfill $\Box$

We now continue our proof using the facing quadruple of strongly separated hyperplanes constructed in the previous lemma. So, up to replacing $\mathfrak{h}_{j}$ by $\mathfrak{h}_{j}'$, we assume that $$(\mathfrak{h}_{j})_{1\le j\le 4}$$ is a facing quadruple of strongly separated halfspaces which moreover do not intersect $\hgotch$.

Exactly as in the proof of the previous lemma, the double skewer lemma implies that there exists $g\in G$ such that $g(\mathfrak{h}_{2}^{\ast})\subsetneq \mathfrak{h}_{1}$ and $h\in G$ such that $\mathfrak{h}_{3} \supsetneq h(\mathfrak{h}_{4}^{\ast})$. Define $U=\mathfrak{h}_{1}\cup \mathfrak{h}_{2}$ and $V=\mathfrak{h}_{3}\cup \mathfrak{h}_{4}$. We now see that we have a Schottky pair:  since $g(V)\subset \mathfrak{h}_{1} \subset \mathfrak{h}_{2}^{\ast}$, for any positive integer $n$, we have $g^{n}(V)\subset \mathfrak{h}_{1}\subset U$. Also, since $g^{-1}(V)\subset \mathfrak{h}_{2}\subset \mathfrak{h}_{1}^{\ast}$ we have $g^{-n}(V)\subset \mathfrak{h}_{2}\subset U$ for any positive integer $n$. Similarly $h^{n}(U)\subset V$ for any non zero integer $n$. The ping-pong lemma implies that $g$ and $h$ generate a free subgroup of $G$. Note that this argument is borrowed from the proof of Theorem F in~\cite{cs}.

Now we observe that $\hgotch \subset U^{c}\cap V^{c}$. If we apply one of the four elements $\{g,g^{-1},h,h^{-1}\}$ to $\hat{\mathfrak{h}}$ we obtain a subset of $U$ or a subset of $V$. This implies that for any nontrivial element $\gamma$ of $\langle g, h\rangle$ we have $\gamma (\hat{\mathfrak{h}}) \subset U$ or $\gamma (\hat{\mathfrak{h}}) \subset V$. In particular $\gamma (\hat{\mathfrak{h}})\cap \hat{\mathfrak{h}}=\emptyset$, and the intersection of the groups $\langle g, h\rangle $ and $Stab_{G}(\hat{\mathfrak{h}})$ is trivial. 

We need to prove something slightly stronger. Namely we are looking for a free subgroup of $G$ which intersects trivially every conjugate of $Stab_{G}(\hgotch)$. So we first make the following observation. Let $x$ be an element of $\langle g,h\rangle$ which is not conjugate to a power of $g$ or of $h$. We will prove that $x$ does not belong to any conjugate of $Stab_{G}(\hgotch)$. Up to changing $x$ into $x^{-1}$ and up to conjugacy, we can assume that this element has the form:
$$x=g^{a_{1}}h^{b_{1}}\cdots g^{a_{r}}h^{b_{r}}$$
with $a_{i}$, $b_{i}$ nonzero integers and $r\ge 1$. This implies that any positive power $x^{m}$ of $x$ satisfies $x^{m}(U)\subset U$. Better, since $x(U)\subset g^{a_{1}}(V)$ it follows from the properties of $g$ and of the $\mathfrak{h}_{i}$'s that the distance of $x(U)$ to the boundary of $U$ is bounded below by some positive number $\delta$. This implies that 
\begin{equation}\label{distu}
d(x^{m}(U),\partial U)\ge m\delta.
\end{equation}
Similarly one proves that 
\begin{equation}\label{distv}
d(x^{-m}(V),\partial V)\ge m\delta'
\end{equation} for any $m\ge 1$ and for some positive number $\delta'$. Suppose now by contradiction that $x$ stabilizes a hyperplane $\hat{\mathfrak{u}}$. Assume that $\hat{\mathfrak{u}}$ does not intersect $U$. For $p\in \hat{\mathfrak{u}}$ we have
$$d(p,U)=d(x^{m}(p),x^{m}(U))$$
for any integer $m\ge 1$. But since $x^{m}(p)$ is never in $U$ the quantity $d(x^{m}(p),x^{m}(U))$ is greater or equal to the distance from $x^{m}(U)$ to $\partial U$. Hence we obtain, using Equation~\eqref{distu}:
$$d(p,U)\ge m\delta$$
which is a contradiction since the left hand side does not depend on $m$. This proves that $\hat{\mathfrak{u}}$ must intersect $U$. In a similar way, using Equation~\eqref{distv}, one proves that $\hat{\mathfrak{u}}$ must intersect $V$. But this contradicts the fact that the halfspaces $(\mathfrak{h}_{j})_{1\le j \le 4}$ are pairwise strongly separated. Hence $x$ does not preserve any hyperplane. In particular $x$ is not contained in any $G$-conjugate of $Stab_{G}(\hgotch)$. 

Now one can consider a normal subgroup $N$ of the free group $\langle g,h\rangle$ which does not contain any nontrivial power of $g$ or $h$ (for instance its derived subgroup). Any finitely generated subgroup $F<N$ has the property that it intersects trivially every $G$-conjugate of $Stab_{G}(\hgotch)$. This proves our claim.

We finally construct two halfspaces $\mathfrak{k}$ and $\mathfrak{l}$ as in the statement of the Theorem. We simply have to repeat the arguments used in the proof of Lemma~\ref{interm}. Exactly as in the proof of that lemma, one can find two halfspaces $\mathfrak{b}_{1}\subset \mathfrak{b}_{2}$ which are contained in $\mathfrak{h}_{j}^{\ast}$ for some $j$ and such that $\hat{\mathfrak{b}}_{1}$ and $\hat{\mathfrak{b}}_{2}$ are strongly separated. We assume that $j=2$ to simplify. We continue as in the proof of the lemma. Applying the element $g$ to $\mathfrak{b}_{1}$ and $\mathfrak{b}_{2}$, we obtain two halfspaces which are contained in $\mathfrak{h}_{1}$: $g(\mathfrak{b}_{1})\subset g(\mathfrak{b}_{2})\subset \mathfrak{h}_{1}$. Since $\mathfrak{h}_{1}\subset \mathfrak{h}_{4}^{\ast}$ and since $h(\mathfrak{h}_{4}^{\ast})\subset \mathfrak{h}_{3}$, one has  
$$hg(\mathfrak{a}_{1})\subset hg(\mathfrak{a}_{2})\subset \mathfrak{h}_{3}.$$
Finally by a similar argument we have: 
$$g^{-1}hg(\mathfrak{b}_{1})\subset g^{-1}hg(\mathfrak{b}_{2})\subset \mathfrak{h}_{2}.$$ We now define $\mathfrak{k}=g(\mathfrak{b}_{1})$ and $\mathfrak{l}=g^{-1}hg(\mathfrak{b}_{1})$. One now checks that $\hgotch$, $\lgotch$, $\kgotch$ are strongly separated exactly as in the end of the proof of Lemma~\ref{interm}.\hfill $\Box$

In what follows, we will use the following definition.

\begin{defi} Let $G$ be a finitely generated group acting on a $\cat$ cubical complex $X$. We say that a hyperplane $\hat{\mathfrak{h}}$ of $X$ is stable for $G$ if the Schreier graph $Stab_{G}(\hat{\mathfrak{h}})\backslash G$ is nonamenable, i.e. satisfies a linear isoperimetric inequality. We say that a halfspace $\mathfrak{h}$ is stable for $G$ if the corresponding hyperplane $\hat{\mathfrak{h}}$ is stable for $G$. 
\end{defi}

Note that to define the Schreier graph of $Stab_{G}(\hat{\mathfrak{h}})\backslash G$, one needs to pick a finite generating set for $G$, but the (non)amenability of this graph is independent of this choice. We refer the reader to~\cite{ikapovich} for the discussion of various equivalent notions of amenability for Schreier graphs. Here we will only need the following:

\begin{lemma}\label{joli} Let $G$ be a group, $H<G$ a subgroup, and $F<G$ a finitely generated free group such that $F\cap gHg^{-1}=\{1\}$ for all $g\in G$. Then any Schreier graph of $H\backslash G$ satisfies a linear isoperimetric inequality i.e. is nonamenable. 
\end{lemma}
{\it Proof.} It is well-known that the Schreier graph of $H\backslash G$ satisfies a linear isoperimetric inequality if and only if the $G$-action on $\ell^{2}(G/H)$ does not have almost invariant vectors. So we will prove this last fact. For this it is enough to prove that the $F$-action on $\ell^{2}(G/H)$ does not have almost invariant vectors. But the hypothesis on $F$ implies that $F$ acts freely on $G/H$. Hence, as an $F$-representation, $\ell^{2}(G/H)$ is isomorphic to a direct sum of copies of the regular representation of $F$ on $\ell^{2}(F)$. This proves the claim.\hfill $\Box$

We will need to record the following corollary of the previous theorem.

\begin{coro}\label{allstable} Under the hypothesis of Theorem~\ref{groupeslibres}, we have:
\begin{itemize}
\item Any halfspace which is part of a facing triple of halfspaces is stable. 
\item For any halfspace $\mathfrak{h}$ which is part of a facing triple of halfspaces, and any finite index subgroup $G_{2}$ of $G$, there exists $\gamma \in G_{2}$ such that $\mathfrak{h}$ and $\gamma (\mathfrak{h})$ are strongly separated. 
\end{itemize}
\end{coro}

\noindent {\it Proof.} This is contained in the proof of the previous theorem. Indeed for the first point of the corollary, we observe that we started the proof of the previous theorem with any facing triple of halfspaces and proved that a given halfspace among the three is stable, as a consequence of Lemma~\ref{joli}. 


For the second point, we consider the triple $\mathfrak{h}, \mathfrak{l}, \mathfrak{k}$ constructed in the previous theorem. We have $\mathfrak{l} \subset \mathfrak{h}^{\ast}$. Applying the double skewer lemma to this last pair, we find $x\in G$ such that $x(\mathfrak{h}^{\ast})\subsetneq \mathfrak{l}$. This implies $x^{n}(\mathfrak{h}^{\ast})\subset \mathfrak{l}$ for all $n\ge 1$. We pick $n_{0}\ge 1$ such that $x^{n_{0}}\in G_{2}$. In particular the hyperplane $x^{n_{0}}(\hgotch)$ is contained in $\mathfrak{l}$. Any hyperplane meeting both $x^{n_{0}}(\hgotch)$ and $\hgotch$ would have to meet $\lgotch$, which is impossible since $\hgotch$ and $\lgotch$ are strongly separated. Hence $x^{n_{0}}(\hgotch)$ and $\hgotch$ are strongly separated.\hfill $\Box$

We will also use the following classical fact. 

\begin{lemma}\label{inisop} Let $M$ be a closed Riemannian manifold with fundamental group $G$. Let $H<G$ be a subgroup and let $M_{1}\to M$ be the covering space associated to $H$. Then the Schreier graph $H\backslash G$ satisfies a linear isoperimetric inequality if and only if $M_{1}$ satisfies a linear isoperimetric inequality.  
\end{lemma}

A proof can be found in~\cite[Ch. 6]{gromovmetric}. Note that the proof in~\cite{gromovmetric} is given only is the case when $H$ is trivial, but the arguments apply in general. Combining Theorem~\ref{groupeslibres}, Corollary~\ref{allstable}, and Lemma~\ref{inisop}, we obtain:

\begin{prop}\label{propstabi} Let $M$ be a closed Riemannian manifold with fundamental group $G$. Suppose that $G$ acts on a $\cat$ cubical complex $X$, satisfying the hypothesis of Theorem~\ref{groupeslibres}. Let $\mathfrak{h}$ be a halfspace of $X$ which is part of a facing triple of halfspaces and let $\hgotch$ be the corresponding hyperplane. Let $M_{\hgotch}$ be the covering space of $M$ corresponding to the subgroup 
$$Stab_{G}(\hgotch)<G.$$
Then, $M_{\hgotch}$ satisfies a linear isoperimetric inequality. 
\end{prop}


\section{Fibering K\"ahler groups acting on $\cat$ cubical complexes}\label{rems}

In this section we first give a criterion to produce fibrations of certain open K\"ahler manifolds over Riemann surfaces (Proposition~\ref{filteredinuse}). Although this criterion is well-known to experts, we explain how to deduce its proof from known results about {\it filtered ends} of K\"ahler manifolds. This serves as a pretext to survey this notion and its applications in the study of K\"ahler groups. Later, in sections~\ref{first} and~\ref{sphm}, we explain how to construct pluriharmonic or plurisubharmonic functions on certain covering spaces of a compact K\"ahler manifold, starting from an action of its fundamental group on a $\cat$ cubical complex. We finally prove Theorem~\ref{fac-pv} in section~\ref{pf}. 

\subsection{Filtered ends}\label{filterede}

The aim of this subsection is to recall the proof of the following classical proposition. 

\begin{prop}\label{filteredinuse} let $M$ be a closed K\"ahler manifold. Let $M_{1}\to M$ be a covering space of $M$ and let $\pi : \widetilde{M}\to M_{1}$ be the universal cover. We assume that there exists a proper, pluriharmonic map $e : M_{1} \to I$ where $I$ is an open interval of $\R$. Let $\widetilde{e}:= e\circ \pi$ be the lift of $e$ to the universal cover $\widetilde{M}$. If some level set of $\widetilde{e}$ is not connected, then there exists a proper holomorphic map with connected fibers 
$$M_{1}\to \Sigma$$
where $\Sigma$ is a Riemann surface. This applies in particular if some levet set of $e$ is not connected. 
\end{prop}

Before turning to the proof of the proposition, we recall briefly various notions of ends in topology and group theory. 

Let $Y$ be a noncompact manifold. Recall that the number of ends of $Y$, denoted by $e(Y)$, is the supremum of the number of unbounded connected components of $Y-K$ where $K$ runs over the relatively compact open sets of $Y$ with smooth boundary. Now if $M$ is a closed manifold, $M_{1}\to M$ a covering space, and $\pi : \widetilde{M}\to M_{1}$ the universal covering space, one can look at the number of ends of various spaces, each of which also admits a purely group theoretical description. 
\begin{itemize}
\item There is the number of ends of $\widetilde{M}$; this is also the number of ends of the fundamental group of $M$.
\item There is the number of ends of the space $M_{1}$, which is an invariant of the pair $\pi_{1}(M_{1})< \pi_{1}(M)$~\cite{houghton,scott}. When this number is greater than $1$, one says that $\pi_{1}(M_{1})$ has codimension $1$ in $\pi_{1}(M)$ or that it is multi-ended. 
\item There is also a third, less known notion, of filtered ends, which can be associated to a continuous map between two manifolds. Here we will only consider the case where this map is the universal covering map. This leads to the following definition.
\end{itemize}

\begin{defi} Let $M_{1}$ be an open manifold and let $\pi : \widetilde{M}\to M_{1}$ be its universal covering space. A filtered end of $M_{1}$ (or of $\pi : \widetilde{M} \to M_{1}$) is an element of the set 
$$\underset{\leftarrow}{{\rm lim}}\, \pi_{0}(\widetilde{M}-\pi^{-1}(K))$$
where $K$ runs over the relatively compact open sets of $M_{1}$ with smooth boundary. The number of filtered ends of $M_{1}$ is denoted by $\widetilde{e}(M_{1})$. 
\end{defi}

As for the previous notions, one can show that in the case where $M_{1}$ covers a compact manifold $M$, this number only depends on the pair $\pi_{1}(M_{1})< \pi_{1}(M)$. Also, one always has $\widetilde{e}(M_{1})\ge e(M_{1})$. The interest in this notion is that the number of filtered ends of $\pi : \widetilde{M}\to M_{1}$ can be greater than $1$ even if $M_{1}$ is $1$-ended. A simple example of this situation is obtained as follows. Take $M$ to be a genus $2$ closed surface and $\Sigma$ to be a subsurface of genus $1$ with one boundary component. The covering space $M_{1}$ of $M$ defined by the subgroup $\pi_{1}(\Sigma)<\pi_{1}(M)$ has this property.

This notion was first introduced by Kropholler and Roller~\cite{kr} in a purely group theoretical context. A topological approach to it was later given in the book by Geoghegan~\cite{ge}. Filtered ends were studied in the context of K\"ahler groups by Gromov and the first author~\cite{delzantgromov}. There, the name {\it cut} was used instead of filtered ends, or rather the term cut was used to indicate the presence of at least two filtered ends for a certain map or covering space, whereas the term {\it Schreier cut} referred to the classical notion of relative ends of a pair of groups in~\cite{delzantgromov}. 

With the notion of filtered ends at our disposal, Proposition~\ref{filteredinuse} will be a simple application of the following theorem.

\begin{theorem}\label{dgnr} Let $M$ be a closed K\"ahler manifold and let $M_{1}\to M$ be an infinite covering space of $M$. If the number of filtered ends of $M_{1}$ is greater or equal to $3$, then there exists a proper holomorphic mapping with connected fibers $M_{1}\to \Sigma$, where $\Sigma$ is a Riemann surface. 
\end{theorem}

This result was proved in~\cite{delzantgromov} under certain additional hypothesis of ``stability". A more general version was later proved by Napier and Ramachandran~\cite{nr2}. Their version includes the theorem stated above but also more general ones which apply to K\"ahler manifolds which are not necessarily covering spaces of a compact manifold. 

Before turning to the proof of Proposition~\ref{filteredinuse}, we make the following easy observation.

\begin{lemma}\label{lieucritnd} Let $V$ be a complex manifold and $f : V \to \R$ be a nonconstant smooth pluriharmonic function. Denote by $Crit(f)$ the set of critical points of $f$. Then for each $t\in \R$, the set $Crit(f)\cap f^{-1}(t)$ is nowhere dense in $f^{-1}(t)$. 
\end{lemma}

The proof is straightforward, once one remembers that the function $f$ is locally the real part of a holomorphic function $F$ and that the critical set of $f$ locally coincides with that of $F$.

\noindent {\it Proof of Proposition~\ref{filteredinuse}.} Note that if $M_{1}$ has at least three ends, the result follows from much older results, see~\cite{nr0}. Thus, we could assume that $M_{1}$ has only two ends, although this is not necessary. 

Let $t_{1}$ be a real number such that $\widetilde{e}^{-1}(t_{1})$ is not connected. Let $x$ and $y$ be two points in distinct connected components of $\widetilde{e}^{-1}(t_{1})$. By Lemma~\ref{lieucritnd}, we can assume that $x$ and $y$ are not critical point of $f$. We claim that at least one of the two sets $$\{ \widetilde{e}>t_{1}\}, \;\;\;\; \{ \widetilde{e} < t_{1}\}$$
is not connected. Let us assume that this is false. Then one can find a path $\alpha$ from $x$ to $y$ contained in the set $\{ \widetilde{e}\ge t_{1}\}$ and a path $\beta$ from $y$ to $x$ contained in the set $\{\widetilde{e} \le t_{1}\}$. We can assume that $\alpha$ and $\beta$ intersect $\widetilde{e}^{-1}(t_{1})$ only at their endpoints. Let $\gamma$ be the path obtained by concatenating $\alpha$ and $\beta$. Let $\gamma_{1} : S^{1}\to M_{1}$ be a smooth map freely homotopic to $\gamma$, intersecting $\widetilde{e}^{-1}(t_{1})$ only at the points $x$ and $y$ and transverse to $\widetilde{e}^{-1}(t_{1})$ there. Since $\widetilde{M}$ is simply connected, we can find a smooth map
$$u : D^{2}\to \widetilde{M}$$
coinciding with $\gamma_{1}$ on the boundary. Pick $\varepsilon >0$ small enough, such that $t_{1}+\varepsilon$  is a regular value of $\widetilde{e}\circ u$. The intersection of the level $\widetilde{e}\circ u = t_{1}+\varepsilon$ with the boundary of the disc is made of two points $a_{\varepsilon}$ and $b_{\varepsilon}$ such that $u(a_{\varepsilon})\to x$ and $u(b_{\varepsilon})\to y$ when $\varepsilon$ goes to $0$. Let $I_{\varepsilon}$ be the connected component of the set $\{\widetilde{e}\circ u = t_{1}+\varepsilon\}$ whose boundary is $\{a_{\varepsilon},b_{\varepsilon}\}$. Up to taking a subsequence, we can assume that $I_{\varepsilon}$ converges as $\varepsilon$ goes to $0$ to a compact connected subset $I_{0}\subset \{\widetilde{e}\circ u =t_{1}\}$. Now the set $u(I_{0})$ is connected, contained in $\widetilde{e}^{-1}(t_{1})$ and contains $x$ and $y$. This contradicts the fact that $x$ and $y$ are in distinct components of $\widetilde{e}^{-1}(t_{1})$. Hence at least one of the two sets $\{ \widetilde{e}>t_{1}\}$ and $\{ \widetilde{e}<t_{1}\}$ is not connected. Note that the maximum principle implies that no connected component of these open sets can be at bounded distance from the set $\widetilde{e}^{-1}(t_{1})$. This is easily seen to imply that $M_{1}$ has at least three filtered ends. The conclusion now follows from Theorem~\ref{dgnr}.\hfill $\Box$

\begin{rem} Note that our Proposition~\ref{filteredinuse} can be applied to the case where $e$ is the primitive of the lift of a harmonic $1$-form $\alpha$ on $M$ with integral periods and where $M_{1}$ is the covering space associated to the kernel of the homomorphism $\pi_{1}(M)\to \Z$ induced by $\alpha$. In this case one recovers a particular case of a result due to Simpson~\cite{simpson}. Hence Theorem~\ref{dgnr} about filtered ends implies this particular case of the result of Simpson. Proposition~\ref{filteredinuse} can also be thought of as a nonequivariant version of Simpson's result. 
\end{rem}

\subsection{The first pluriharmonic map}\label{first}

We now consider a ${\rm CAT}(0)$ cubical complex $X$, assumed to be irreducible and locally finite. Let $M$ be a compact K\"ahler manifold and let $\Gamma$ be its fundamental group. We consider a homomorphism
$$\varrho : \Gamma \to {\rm Aut}(X).$$
We suppose that the $\Gamma$-action on $X$ satisfies the three hypothesis of Theorem~\ref{fac-pv}, although this will not be used untill section~\ref{sphm}. We denote by $\widetilde{M}$ the universal cover of $M$. If $\hgotch$ is a hyperplane of $X$ we will write $M_{\hgotch}$ for the quotient of $\widetilde{M}$ by the group $Stab_{\Gamma}(\hgotch)$:
\begin{equation}\label{notationrev}
M_{\hgotch}=\widetilde{M}/Stab_{\Gamma}(\hgotch).
\end{equation}
Since $X$ is contractible, one can choose a Lipschitz $\varrho$-equivariant map $\psi : \widetilde{M}\to X$. We fix a hyperplane $\hat{\mathfrak{h}}$ of $X$ and a halfspace $\mathfrak{h}$ associated to $\hat{\mathfrak{h}}$. Define a function $d_{\mathfrak{h}}$ on $X$ by 
$$d_{\mathfrak{h}}(x)=\left\{ \begin{array}{ll}
d(x,\hat{\mathfrak{h}}) & if\; x\in \mathfrak{h}, \\
-d(x,\hat{\mathfrak{h}}) & if \; x\in \mathfrak{h}^{\ast}. \\
\end{array}\right.$$ 
Let $f : \widetilde{M}\to \R$ be the composition $f=d_{\mathfrak{h}}\circ \psi$. Finally let $\overline{f}$ be the function induced by $f$ on $M_{\hgotch}$.

\begin{prop}\label{proprete} The map $\overline{f}$ is proper. In particular the manifold $M_{\hgotch}$ has at least two ends. 
\end{prop}
{\it Proof.} In the following proof, we let $H=Stab_{G}(\hgotch)$. We pick a point $x_{0}$ in $\widetilde{M}$. We consider the map
$$F : H\backslash G \to \R$$
defined by $F(Hg)=f(g\cdot x_{0})$. It is enough to prove that $F$ is proper, this is easily seen to imply that $\overline{f}$ is proper. So let $(g_{n})$ be a sequence of elements of $G$ such that 
\begin{equation}\label{borne}
\vert F(Hg_{n})\vert \le C,\end{equation}
for some constant $C$. We need to prove that the set $\{Hg_{n}\}$ is a finite subset of $H\backslash G$. But equation~\eqref{borne} implies that the distance of $\varrho (g_{n})(\psi(x_{0}))$ to $\hat{\mathfrak{h}}$ is bounded by $C$. Put  differently this says that each of the hyperplanes $(\varrho (g_{n})^{-1}(\hat{\mathfrak{h}}))$ intersects the ball of radius $C$ centered at $\psi (x_{0})$. By our hypothesis of local finiteness of $X$ this implies that the family of hyperplanes $(\varrho (g_{n})^{-1}(\hgotch))$ is {\it finite}. Hence there exists element $x_{1}, \ldots , x_{r}$ in $G$ such that for all $n$ there exists $i$ such that $\varrho (g_{n})^{-1} (\hgotch ) = \varrho (x_{j})(\hgotch)$. This implies $Hg_{n}=Hx_{i}^{-1}$.\hfill $\Box$

In what follows we will denote by $Ends(V)$ the space of ends of an open manifold $V$. We recall that $V\cup Ends(V)$ carries a natural topology, see for instance~\cite[Ch.6]{dk}. We now make use of the following classical theorem. 

\begin{theorem}\label{tresjoli} Let $V$ be a complete Riemannian manifold of bounded geometry and satisfying a linear isoperimetric inequality. Let $\chi : Ends(V)\to \{-1,1\}$ be a continuous map. Then there exists a unique continuous map $h_{\chi} : V\cup Ends(V)\to [-1,1]$ which extends $\chi$, is harmonic on $V$, and has finite energy:
$$\int_{V}\vert \nabla h_{\chi}\vert^{2}<\infty.$$
\end{theorem}

Note that we will not need the precise definition of a Riemannian manifold of bounded geometry; the only important point we need is that a covering space of a closed Riemannian manifold is of bounded geometry. This theorem was proved by Kaimanovich and Woess~\cite{kw} and by Li and Tam~\cite{litam} independently. A simple beautiful proof due to M. Ramachandran can be found in the paper by Kapovich~\cite[\S 9]{mkapovich}. In the case when $V$ is K\"ahler, the map $h_{\chi}$ is pluriharmonic. This follows from a standard integration by part argument, which is valid here because $h_{\chi}$ has finite energy, see Lemma 3.1 in~\cite{li} or ~\cite[\S 1.1.B]{khgromov}.

We now assume that $\hgotch$ is stable for $\Gamma$ and apply the above to the manifold $V=M_{\hgotch}$. The proper function $\overline{f}$ defines a partition of the space of ends of this manifold into two open sets: the set of ends of $\{\overline{f}\ge 0\}$ and the set of ends of $\{\overline{f}\le 0\}$. Let $\chi : Ends(V)\to \{-1,1\}$ be the function taking the value $1$ on the first open set and $-1$ on the second open set. Since $M_{\hgotch}$ satisfies a linear isoperimetric inequality, we obtain:
 
\begin{prop}\label{tropcool} There exists a unique function
$$u_{\mathfrak{h}} : M_{\hgotch} \to (-1,1)$$
which is pluriharmonic, of finite energy and satisfies $u_{\mathfrak{h}}(x)\to 1$ when $x\to \infty$ and $\overline{f}(x) >0$ and $u_{\mathfrak{h}}(x)\to -1$ when $x\to \infty$ and $\overline{f}(x) <0$.
\end{prop}

We denote by $p$ the projection from $\widetilde{M}$ to $M_{\hgotch}$. For $t_{0}\in (0,1)$ we define: $$U_{\mathfrak{h},t_{0}} =\{x\in \widetilde{M}, u_{\mathfrak{h}}(p(x))>t_{0}\}.$$ We have: 

\begin{lemma}\label{bhd} 
If $t_{0}$ is close enough to $1$ the set $U_{\mathfrak{h},t_{0}}$ satisfies $\psi  (U_{\mathfrak{h},t_{0}}) \subset \mathfrak{h}$.
\end{lemma}
{\it Proof.} By the continuity of $u_{\mathfrak{h}}$ on $M_{\hgotch}\cup Ends(M_{\hgotch})$, we have that if $t_{0}$ is close enough to $1$, one has the inclusion $\{u_{\mathfrak{h}}> t_{0}\}\subset \{\overline{f}\ge 0\}$. This implies the conclusion of the lemma.\hfill $\Box$

We will ultimately produce a fibration of the manifold $M_{\hgotch}$ onto a Riemann surface. To achieve this, we will construct a second plurisubharmonic function defined on $M_{\hgotch}$, which is not a function of $u_{\mathfrak{h}}$. We build this second map in the next section.

\subsection{A second plurisubharmonic map}\label{sphm}

We consider a facing triple of strongly separated halfspaces $\mathfrak{h}$, $\mathfrak{k}$, $\mathfrak{l}$ as given by Theorem~\ref{groupeslibres}. According to Corollary~\ref{allstable},  $\mathfrak{h}$, $\mathfrak{k}$ and $\mathfrak{l}$ are stable for $\Gamma$. Thus, we can apply the results of section~\ref{first} to anyone of the three covering spaces $M_{\hgotch}$, $M_{\kgotch}$, $M_{\lgotch}$. We recall that these covering spaces were defined by Equation~\eqref{notationrev}. We also make use of the various lemmas proved in section~\ref{essh}. In particular we consider the subgroup 
$$A< Stab_{\Gamma}(\hgotch)$$
generated by $\Sigma=\{h\in Stab_{\Gamma}(\hgotch), h(\mathfrak{k})\cap \mathfrak{k}\neq \emptyset \}$. If $a\in A$, the hyperplane $a(\kgotch)$ is stable for $\Gamma$, as $\kgotch$ is. Thus, we can also apply the results of section~\ref{first} to the covering space $M_{a(\kgotch)}$. So we denote by 
$$u_{\mathfrak{h}} : M_{\hgotch}\to (-1,1)\;\;\;\;\; {\rm and} \;\;\;\;\; u_{a(\mathfrak{k})} : M_{a(\kgotch)} \to (-1,1)$$
the proper pluriharmonic maps of finite energy provided by Proposition~\ref{tropcool}.

For each $a\in A$, let $\widetilde{u_{a(\mathfrak{k})}}$ be the lift of $u_{a(\mathfrak{k})}$ to $\widetilde{M}$. Now define a function $\beta_{\mathfrak{k}} : \widetilde{M}\to \R$ by:
$$\beta_{\mathfrak{k}}(x)=\underset{a\in A}{{\rm sup}} \, \widetilde{u_{a(\mathfrak{k})}}(x).$$
Proposition~\ref{finitude} implies that the $A$-orbit of $\mathfrak{k}$ is finite, hence the supremum above is actually the supremum of a finite number of smooth pluriharmonic functions. Hence $\beta_{\mathfrak{k}}$ is a continuous plurisubharmonic function. We now fix some constant $C$ in $(0,1)$. Let $V_{\mathfrak{k},C}=\{\beta_{\mathfrak{k}}>C\}$.

\begin{lemma} \label{lemmeinterm}
\begin{enumerate}
\item The function $\beta_{\mathfrak{k}} : \widetilde{M}\to \R$ is $A$-invariant. In particular the open set $V_{\mathfrak{k},C}\subset \widetilde{M}$ is $A$-invariant.  
\item For each $a\in A$, there exists $C_{a}\in(0,1)$ such that $\psi (\{\widetilde{u_{a(\mathfrak{k})}}>C_{a}\})\subset a(\mathfrak{k})$.
\item If $C$ is close enough to $1$, one has $\psi(V_{\mathfrak{k},C})\subset \cup_{a\in A} a(\mathfrak{k})$.\label{lemmeinterm:3}
\end{enumerate}
\end{lemma}
{\it Proof.} Note that by the uniqueness of the map $u_{a(\mathfrak{k})}$, we have:
$$\widetilde{u_{a(\mathfrak{k})}}=\widetilde{u_{\mathfrak{k}}}\circ a^{-1}.$$
So the function $\beta_{\mathfrak{k}}$ is also equal to
$$\underset{a\in A}{{\rm sup}}\, \widetilde{u_{\mathfrak{k}}}\circ a^{-1},$$
which is $A$-invariant. This proves the first point. The second point was already proved in Lemma~\ref{bhd} for the harmonic function $u_{\mathfrak{h}} : M_{\hgotch}\to (-1,1)$ and the proof is similar here. As for the third point, if $A_{1}\subset A$ is a finite set such that $\{a(\mathfrak{k})\}_{a\in A}=\{a(\mathfrak{k})\}_{a\in A_{1}}$, the constant $C=\underset{a\in A_{1}}{{\rm max}} \, C_{a}$ has the desired property.\hfill $\Box$

We now assume that the constant $C$ satisfies the conclusion of the previous lemma. Note that the open set $V_{\mathfrak{k},C}$ naturally defines an open set $V^{\ast}_{\mathfrak{k},C}=V_{\mathfrak{k},C}/A$ inside $\widetilde{M}/A$. This open set actually embeds into the quotient of $\widetilde{M}$ by the bigger subgroup $Stab_{\Gamma}(\hgotch)$. Namely:

\begin{prop} The natural map $\pi : \widetilde{M}/A\to M_{\hgotch}$ is injective on the closure of $V^{\ast}_{\mathfrak{k},C}$. 
\end{prop}
{\it Proof.} We will check that if $h\in Stab_{\Gamma}(\hgotch)-A$ then $h(\overline{V_{\mathfrak{k},C}})\cap \overline{V_{\mathfrak{k},C}}$ is empty. Let us suppose by contradiction that this is not empty for some $h\in Stab_{\Gamma}(\hgotch)-A$. Let $x\in \overline{V_{\mathfrak{k},C}}$ be such that $h(x)\in \overline{V_{\mathfrak{k},C}}$. By Lemma~\ref{lemmeinterm}~\eqref{lemmeinterm:3}, there exist $a_{1}, a_{2}\in A$ such that $\psi(x)\in a_{1}(\mathfrak{k})$ and $\psi(h(x))\in a_{2}(\mathfrak{k})$. This implies that $h(a_{1}(\mathfrak{k}))\cap a_{2}(\mathfrak{k})\neq \emptyset$. By Lemma~\ref{separea}, $h$ must be in $A$, a contradiction.\hfill $\Box$

\medskip

We now define a function $\gamma : M_{\hgotch} \to \R_{+}$ as follows. If $x\in \pi(V^{\ast}_{\mathfrak{k},C})$, define $\gamma (x)=\beta_{\mathfrak{k}}(y)-C$, where $y$ is a lift of $x$ inside $V_{\mathfrak{k},C}$. If $x\notin \pi (V^{\ast}_{\mathfrak{k},C})$, define $\gamma (x)=0$. We have:

\begin{lemma}
If $C$ is close enough to $1$, the function $\gamma$ is continuous and plurisubharmonic. Hence, $V^{\ast}_{\mathfrak{k},C}$ is a plurimassive set in the sense of~\cite[Def. 1.1]{nr2}.
\end{lemma}
{\it Proof.} As in the previous Proposition, we denote by $\pi$ the projection $\widetilde{M}/A\to M_{\hgotch}$. Up to taking $C$ closer to $1$, we can always assume that there exists $C_{0}\in (0,C)$ such that the conclusion of the third point of Lemma~\ref{lemmeinterm} holds for $C_{0}$. We claim that $F:=\pi(\overline{V_{\mathfrak{k},C}}/A)$ is closed. We first conclude the proof of the lemma using this claim. We consider the two open sets
$$\pi(V_{\mathfrak{k},C_{0}}/A)\;\;\;\;\; and \;\;\;\;\; M_{\hgotch}-F.$$
They cover $M_{\hgotch}$, hence it is enough to prove that $\gamma$ is continuous and plurisubharmonic on each of them. On $M_{\hgotch}-F$, $\gamma$ is $0$ so there is nothing to prove. On $\pi(V_{\mathfrak{k},C_{0}}/A)$, the function $\gamma$ is constructed as follows: one considers the plurisubharmonic function ${\rm max}(\beta_{\mathfrak{k}}-C,0)$ on $\widetilde{M}$. It descends to a function $q : \widetilde{M}/A\to \R_{+}$. The restriction of $\gamma$ to $\pi(V_{\mathfrak{k},C_{0}}/A)$ is obtained by composing the inverse of the map $\pi : V_{\mathfrak{k},C_{0}}/A \to \pi(V_{\mathfrak{k},C_{0}}/A)$ with $q$. Hence it is continuous and plurisubharmonic.

We finally prove that $F:=\pi(\overline{V_{\mathfrak{k},C}}/A)$ is closed. It is enough to see that 
$$\bigcup_{h\in Stab_{\Gamma}(\hgotch)}h(\overline{V_{\mathfrak{k},C}})$$
is closed. For this it is enough to check that there exists $\varepsilon >0$ such that for any $h\in Stab_{\Gamma}(\hgotch)-A$, the distance from $h(\overline{V_{\mathfrak{k},C}})$ to $\overline{V_{\mathfrak{k},C}}$ is greater or equal to $\varepsilon$. Applying the map $\psi$, we see that it is enough to find a positive lower bound, independent of $h\in Stab_{\Gamma}(\hgotch)-A$, for the distance from $h(U)$ to $U$ in $X$, where $U=\cup_{a\in A}a(\mathfrak{k})$, as in section~\ref{essh}. But this follows from the fact that there is a uniform positive lower bound for the distance between two disjoint halfspaces in a $\cat$ cubical complex.\hfill $\Box$

We started this section considering a facing triple of strongly separated hyperplanes $\mathfrak{h}$, $\mathfrak{k}$ and $\mathfrak{l}$. So far, we only used $\mathfrak{h}$ and $\mathfrak{k}$. In the next proposition, we make use of the third hyperplane $\mathfrak{l}$.  

\begin{prop}\label{supercool} Assume that the level sets of $\widetilde{u_{\mathfrak{h}}}$ (the lift of $u_{\mathfrak{h}}$ to $\widetilde{M}$) are connected. Then, there exists a finite cover $M_{2} \to M_{\hgotch}$ (possibly equal to $M_{\hgotch}$) and a continuous plurisubharmonic function $\delta : M_{2}\to \R_{+}$ such that there exists a level set of the lift of $u_{\mathfrak{h}}$ to $M_{2}$ on which $\delta$ is not constant. 
\end{prop}
{\it Proof.} Note that the function $u_{\mathfrak{h}}$ is surjective. Assume that the conclusion of the proposition is false when $M_{2}=M_{\hgotch}$ and $\delta=\gamma$. Then there exists a function
$$\varphi : (-1,1) \to \R_{+}$$
such that $\gamma =\varphi \circ u_{\mathfrak{h}}$. We claim that the function $\varphi$ is continuous, convex, and vanishes on $[t_{0},1)$ for $t_{0}$ close enough to $1$. Let us prove these claims. 

First we note that $\varphi$ is continuous in a neighborhood of every real number $t$ such that $u_{\mathfrak{h}}^{-1}(t)$ is not contained in the critical set of $u_{\mathfrak{h}}$. But every real number $t$ has this property according to Lemma~\ref{lieucritnd}. Hence $\varphi$ is continuous on $(-1,1)$.  Note that this argument proves that for each $t\in (-1,1)$, one can find $q\in M_{\hgotch}$  and local coordinates $(z_{1},\ldots , z_{n})$ centered at $q$ such that $\gamma(z)=\varphi(u_{\mathfrak{h}}(q)+Re(z_{1}))$. The convexity of $\varphi$ then follows from~\cite[5.13]{demailly}. For the last claim, we pick a point $x\in \widetilde{M}$ such that $u_{\mathfrak{h}}(p(x))>t_{0}$. Here and as before $p$ denotes the projection $\widetilde{M}\to M_{\hgotch}$. If $t_{0}$ is close enough to $1$, Lemma~\ref{bhd} implies that $\psi (x)\in \mathfrak{h}$. Since $\mathfrak{h}$ and
$$\bigcup_{a\in A}a(\mathfrak{k})$$
are disjoint, Lemma~\ref{lemmeinterm} implies that $h(x)\notin V_{\mathfrak{k}, C}$ for any $h$ in $Stab_{\Gamma}(\hgotch)$. Hence $\gamma (p(x))=0$. This proves that $\varphi$ vanishes on $[t_{0},1)$.

These three properties of $\varphi$ imply that the level sets of $\varphi$ are connected. In fact the level $\varphi = c$ is a point for $c>0$ and is an interval of the form $[a_{0},1)$ for $c=0$. This implies, together with the hypothesis on the level sets of $\widetilde{u_{\mathfrak{h}}}$, that the level sets of $\gamma \circ p : \widetilde{M} \to \R_{+}$ are connected. But this implies $Stab_{\Gamma}(\hgotch)=A$. Now since the $A$-orbit of $\mathfrak{k}$ is finite, the group
$$H_{2}=Stab_{\Gamma}(\hgotch)\cap Stab_{\Gamma}(\kgotch)$$
is of finite index in $Stab_{\Gamma}(\hgotch)$. Let $M_{2}\to M_{1}$ be the corresponding cover and $u_{2}$ be the lift of $u_{\mathfrak{h}}$ to $M_{2}$. Let $\delta$ be the lift of the function $u_{\mathfrak{k}}$ to $M_{2}$. We claim that $\delta$ satisfies the conclusion of the proposition. 

If this is not the case, then $\delta$ is a function of $u_{2}$. As before one sees that $\delta =\varphi_{0} \circ u_{2}$ where $\varphi_{0}$ is continuous and convex. Actually, since $\delta$ is pluriharmonic and not only plurisubharmonic, $\varphi_{0}$ must be affine. So there exists real numbers $a$ and $b$ such that:
$$\delta = au_{2}+b.$$
Since $u_{2}$ and $\delta$ both take values into $(-1,1)$ and are onto, one sees that $b$ must be $0$ and that $a=\pm 1$. We now obtain a contradiction from this fact, making use of the hyperplane $\mathfrak{l}$. Let $s_{n}$ be a sequence of points of $\widetilde{M}$ such that $d(\psi (s_{n}),\hat{\mathfrak{l}})$ goes to infinity and such that $\psi (s_{n})\in \mathfrak{l}$. Such a sequence exists because the action is essential. Since $d(\psi(s_{n}),\mathfrak{h})\ge d(\psi (s_{n}),\hat{\mathfrak{l}})$ we must have that $\overline{f}(p(s_{n}))\to -\infty$, hence also $u_{\mathfrak{h}}(p(s_{n}))\to -1$. 

In the next paragraph, we denote by $x\mapsto [x]$ the covering map $\widetilde{M}\to M_{2}$. 

If $a=-1$, one sees that $\delta ([s_{n}])\to 1$ which implies that $\psi (s_{n})\in \mathfrak{k}$ for $n$ large enough. Hence $\Psi (s_{n})\in \mathfrak{k}\cap \mathfrak{l}$. This is a contradiction since $\mathfrak{k}$ and $\mathfrak{l}$ are disjoint. If $a=1$ we argue in a similar way with the pair $(\mathfrak{h},\mathfrak{k})$. We take a sequence $q_{n}$ of points of $\widetilde{M}$ such that $\psi (q_{n})\in \mathfrak{h}$ and $d(\psi(q_{n}),\hgotch)\to \infty$. This implies $u_{2}([q_{n}])\to 1$. Since $\delta = u_{2}$ this implies that $\psi (q_{n})\in \mathfrak{k}$ for $n$ large enough. Since $\mathfrak{k}\cap \mathfrak{h}$ is empty we get a contradiction again. This proves the proposition.\hfill $\Box$

\subsection{Producing fibrations}\label{pf}

We continue with the notations and hypothesis from section~\ref{sphm}. Our aim is now to prove the following:

\begin{prop}\label{onemore} The manifold $M_{\hgotch}$ {\it fibers}: there exists a proper holomorphic map $M_{\hgotch}\to \Sigma$ with connected fibers onto a Riemann surface $\Sigma$. 
\end{prop}

\noindent {\it Proof.} By Proposition~\ref{filteredinuse}, we can assume that the level sets of the lift of the map
$$u_{\mathfrak{h}} : M_{\hgotch}\to (-1,1)$$
to the universal cover are connected, otherwise the conclusion is already known. By Proposition~\ref{supercool}, we can first replace $M_{\hgotch}$ by a finite cover $p : M_{2}\to M_{\hgotch}$ such that there exists a function $\delta : M_{2} \to \R_{+}$ which is continuous, plurisubharmonic and not constant on the set $\{u_{\mathfrak{h}}\circ p =t\}$ for some real number $t$. Note that if this is true for some number $t$, one can always find a number $t'$, close to $t$, which has the same property and moreover satisfies that $t'$ is a regular value of $u_{\mathfrak{h}}\circ p$. Let $\mathscr{F}$ be the foliation of $M_{2}$ defined by the $(1,0)$ part of the differential of $u_{\mathfrak{h}}\circ p$. Note that $d(u_{\mathfrak{h}}\circ p)$ might have zeros; we refer to~\cite[p.55]{abckt} for the precise definition of $\mathscr{F}$. We now consider the manifold 
$$X_{t'}=(u_{\mathfrak{h}}\circ p)^{-1}(t').$$
It is invariant by the foliation $\mathscr{F}$. Hence on $X_{t'}$ we have a real codimension $1$ foliation defined by the nonsingular closed $1$-form which is the restriction of $d^{\mathbb{C}}(u_{\mathfrak{h}}\circ p)$ to $X_{t'}$. Such a foliation has all its leaves closed or all its leaves dense; this is an elementary particular case of the theory of Riemannian foliations. We must thus show that the restriction of $\mathscr{F}$ to $X_{t'}$ cannot have all its leaves dense. Let $q$ be a point where the restriction of $\delta$ to $X_{t'}$ reaches its maximum $m$. Let $L(q)$ be the leaf of $\mathscr{F}$ through $q$. The maximum principle implies that $\delta$ is constant on $L(q)$. Hence the closure of $L(q)$ is contained in the set $\delta =m$. Since $\delta$ is not constant on $X_{t'}$, $L(q)$ cannot be dense in $X_{t'}$, hence is closed. We have found a compact leaf of the foliation $\mathscr{F}$. This compact leaf projects to a compact leaf of the foliation $\mathscr{F}_{\hgotch}$ defined by $du_{\mathfrak{h}}^{1,0}$ on $M_{\hgotch}$. But it is now classical that this implies that the foliation $\mathscr{F}_{\hgotch}$ is actually a fibration. See for instance~\cite[(7.4)]{carlsontoledo} or~\cite[\S 4.1]{delzantgromov}.\hfill $\Box$

Now we will apply the following result:

{\it Let $V$ be a closed K\"ahler manifold. Assume that $V$ has a covering space $V_{1}\to V$ which admits a proper holomorphic mapping to a Riemann surface, with connected fibers: $\pi_{1} : V_{1}\to \Sigma_{1}$. Then $\pi_{1}$ descends to a finite cover of $V$: there exists a finite cover $V_{2}$ of $V$ such that $V_{1}\to V$ decomposes as $V_{1}\to V_{2}\to V$, $V_{2}$ admits a holomorphic mapping $\pi_{2} : V_{2} \to \Sigma_{2}$ with connected fibers and there exists a holomorphic map $\Sigma_{1}\to \Sigma_{2}$ which makes the following diagram commutative:
\begin{displaymath}
\xymatrix{V_{1} \ar[r] \ar[d] & \Sigma_{1} \ar[d] \\
V_{2} \ar[r] & \Sigma_{2}. \\}
\end{displaymath} }

This fact is now well-known, we refer the reader to~\cite[\S 5.6]{delzantgromov} or~\cite[Prop. 4.1]{nr1} for a proof. Applying this result to $V=M$ and to the cover $V_{1}=M_{\hgotch}$ we obtain a finite cover $M_{2}\to M$ and a fibration $\pi_{2} : M_{2}\to \Sigma$ onto a Riemann surface. By replacing $M_{2}$ by another finite cover, we can assume that the fundamental group of a smooth fiber of $\pi_{2}$ surjects onto the kernel of the map $(\pi_{2})_{\ast} : \pi_{1}(M_{2})\to \pi_{1}(\Sigma)$. Note that a smooth fiber of the fibration of $M_{\hgotch}$ obtained in Proposition~\ref{onemore} projects onto a smooth fiber of $\pi_{2}$. This implies that the normal subgroup $$N:=Ker((\pi_{2})_{\ast} : \pi_{1}(M_{2})\to \pi_{1}(\Sigma))$$
is contained in the stabilizer of $\hgotch$ inside $\pi_{1}(M_{2})$. In what follows, we write $\Gamma_{2}=\pi_{1}(M_{2})$. To conclude the proof of Theorem~\ref{fac-pv}, we only have to establish the next proposition. 

\begin{prop} The normal subgroup $N$ acts as an elliptic subgroup of ${\rm Aut}(X)$ i.e. the fixed point set of $N$ in $X$ is nonempty. 
\end{prop}
{\it Proof.} We know that $N$ is contained in the group $Stab_{\Gamma_{2}}(\hgotch)$. But since $N$ is normal in $\Gamma_{2}$, $N$ is contained in 
$$Stab_{\Gamma_{2}}(\gamma(\hgotch))$$
for all $\gamma$ in $\Gamma_{2}$. By Corollary~\ref{allstable}, we can pick $\gamma \in \Gamma_{2}$ such that $\hgotch$ and $\gamma (\hgotch)$ are strongly separated. Since $N$ preserves $\hgotch$ and $\gamma (\hgotch)$ it must preserve the projection of $\hgotch$ onto $\gamma(\hgotch)$, which is a point according to Proposition~\ref{projunique}. Hence $N$ is elliptic.\hfill $\Box$

\section{Cubulable K\"ahler manifolds and groups}\label{ckmag}

\subsection{Cubulable K\"ahler groups}

We first recall a few definitions concerning finiteness properties of groups as well as a result by Bridson, Howie, Miller and Short~\cite{bhms2002}. These will be used in the proof of Theorem~\ref{fac-complet}. A group $G$ is of type ${\rm FP}_{n}$ if there is an exact sequence 
$$P_{n}\to P_{n-1}\to \cdots \to P_{0}\to \Z \to 0$$
\noindent of $\Z G$-modules, where the $P_{i}$ are finitely generated and projective and where $\Z$ is considered as a trivial $\Z G$-module. It is of type ${\rm FP}_{\infty}$ if it is of type ${\rm FP}_{n}$ for all $n$. See~\cite[\S VIII.5]{browncoho} for more details on these notions. We simply mention that the fundamental group of a closed aspherical manifold is of type ${\rm FP}_{\infty}$. 

It is proved in~\cite{bhms2002} that if $H_{1}, \ldots , H_{n}$ are either finitely generated free groups or surface groups and if $G$ is a subgroup of the direct product 
$$H_{1}\times \cdots \times H_{n}$$
which is of type ${\rm FP}_{n}$, then $G$ is virtually isomorphic to a direct product of at most $n$ groups, each of which is either a surface group of a finitely generated free group. We also refer the reader to~\cite{bridsonhowie,bhms2009} for more general results. Note that the idea of applying the results from~\cite{bhms2002} to K\"ahler groups is not new. This possibility was already discussed in~\cite{bridsonhowie,bhms2009}, and put into use in~\cite{py}. We also mention here that there exist K\"ahler groups which are subgroups of direct products of surface groups but which are not of type ${\rm FP}_{\infty}$, see~\cite{dps,cli}.

We now prove Theorem~\ref{fac-complet}. So let $\Gamma$ and $X$ be as in the statement of the theorem. Let 
$$X= X_{1}\times \cdots \times X_{r}$$
be the decomposition of $X$ into a product of irreducible $\cat$ cubical complexes. There is a finite index subgroup $\Gamma_{1}$ of $\Gamma$ which preserves this decomposition. Note that the action of $\Gamma_{1}$ on each of the $X_{i}$ is essential since the original action is essential on $X$. We will make use of the following two results from~\cite{cs}:
\begin{itemize}
\item The group $\Gamma_{1}$ contains an element $\gamma_{0}$ which acts as a rank $1$ isometry on each irreducible factor. 
\item The group $\Gamma_{1}$ contains a copy of $\Z^{r}$. 
\end{itemize}
See~\cite{cs}, Theorem C and Corollary D for these statements. We recall here that a rank $1$ isometry of a $\cat$ space is a hyperbolic isometry none of whose axis bounds a flat half-plane.

\begin{prop}\label{tuerparab} Let $i\in \{1, \ldots , r\}$. Exactly one of the following two cases occurs:
\begin{enumerate}
\item The action of $\Gamma_{1}$ preserves a geodesic line in $X_{i}$. 
\item The action of $\Gamma_{1}$ has no invariant Euclidean flat in $X_{i}$ and fixes no point in the viusal boundary of $X_{i}$. 
\end{enumerate}
\end{prop}
{\it Proof.} Since $\Gamma_{1}$ contains a rank $1$ isometry, if $\Gamma_{1}$ preserves a Euclidean flat in $X_{i}$, this flat must be a geodesic line. According to Proposition 7.3 from~\cite{cs} (see also the proof of that proposition), if $\Gamma_{1}$ does not preserve a geodesic line in $X_{i}$, it does not have any fixed point in the visual boundary of $X_{i}$.\hfill $\Box$

Up to replacing $\Gamma_{1}$ by a finite index subgroup, we assume that whenever $\Gamma_{1}$ preserves a geodesic line $L_{i}$ in some $X_{i}$, it acts by translations on it. In this case, the translation group is discrete as follows from~\cite{bridson} for instance. Hence the action of $\Gamma_{1}$ on $L_{i}$ factors through a homomorphism to $\Z$. 

We now continue the proof of Theorem~\ref{fac-complet}. We change the numbering of the factors $X_{i}$'s so that for $1\le i \le k$, $\Gamma_{1}$ preserves a geodesic line in $X_{i}$ and acts by translations on it, whereas for $i>k$, it fixes no point in the visual boundary of $X_{i}$ and preserves no flat. Hence for $j>k$, the $\Gamma_{1}$-action on $X_{j}$ satisfies all the hypothesis of Theorem~\ref{fac-pv}. Using Theorem~\ref{fac-pv} we get a finite index subgroup $\Gamma_{2}<\Gamma_{1}$ and a homomorphism 
\begin{equation}\label{ssttrr}
\psi : \Gamma_{2}\to \Z^{k}\times \pi_{1}(\Sigma_{1})\times \cdots \times \pi_{1}(\Sigma_{r-k})
\end{equation}
with the following properties:
\begin{enumerate}
\item The homomorphism $\Gamma_{2}\to \pi_{1}(\Sigma_{j})$ induced by $\psi$ is surjective for each $k+1\le j \le r$.
 \item For each $i$, there exists a convex cobounded subset $Y_{i}\subset X_{i}$ on which the $\Gamma_{2}$-action factors through the coordinate number $i$ of the homomorphism $\Gamma_{2}\to \Z^{k}$ if $i\le k$ or through the homomorphism $\Gamma_{2} \to \pi_{1}(\Sigma_{i-k})$ if $i>k$. 
\end{enumerate}

\noindent Since the action of $\Gamma_{2}$ is proper, the kernel $N$ of $\psi$ is finite. By replacing $\Gamma_{2}$ by a finite index subgroup $\Gamma_{3}$, we can assume that $N\cap \Gamma_{3}$ is central in $\Gamma_{3}$. We thus have a central extension: 
\begin{displaymath}
\xymatrix{ \{0\} \ar[r] & N\cap \Gamma_{3} \ar[r] & \Gamma_{3} \ar[r] & \psi(\Gamma_{3}) \ar[r] & \{0\}. \\
}
\end{displaymath}


\begin{lemma}\label{fpfpfp} The group $\psi(\Gamma_{3})$ is of type ${\rm FP}_{\infty}$. 
\end{lemma}
{\it Proof.} Let $Y\subset X$ be the fixed point set of $N$. This is a subcomplex of the first cubical subdivision of $X$. Since $\psi(\Gamma_{3})$ is torsion-free, it acts freely on $Y$. The quotient $Y/\psi(\Gamma_{3})$ is a finite complex, hence $\psi(\Gamma_{3})$ is of type ${\rm FL}$, in particular ${\rm FP}_{\infty}$. See~\cite[VIII.6]{browncoho} for the definition of the ${\rm FL}$ condition and its relation to the ${\rm FP}_{n}$ and ${\rm FP}_{\infty}$ conditions.\hfill $\Box$

Now the result of~\cite{bhms2002} implies that $\psi(\Gamma_{3})$ itself is isomorphic to a direct product of surface groups and finitely generated free groups. No non-Abelian free factor can appear, but since this is not needed for the next lemma, we will postpone a little bit the explanation of this fact.   

\begin{lemma} The group $\Gamma_{3}$ has a finite index subgroup $\Gamma_{4}$ which does not intersect $N$. In particular $\Gamma_{4}\simeq \psi(\Gamma_{4})$.
\end{lemma}
{\it Proof.} It is enough to prove that the central extension of $\psi (\Gamma_{3})$ by $N\cap \Gamma_{3}$ appearing above becomes trivial on a finite index subgroup of $\psi (\Gamma_{3})$. Being isomorphic to a direct product of surface groups and free groups, $\psi(\Gamma_{3})$ has torsion-free $H_{1}$. Hence the universal coefficient theorem implies that
$$H^{2}(\psi (\Gamma_{3}),N\cap \Gamma_{3})$$
is isomorphic to $Hom(H_{2}(\psi(\Gamma_{3}),\Z),N\cap \Gamma_{3})$. If $H<\psi(\Gamma_{3})$ is a subgroup of finite index such that every element in the image of the map $H_{2}(H,\Z)\to H_{2}(\psi (\Gamma_{3}),\Z)$ is divisible by a large enough integer $p$, the pull-back map $Hom(H_{2}(\psi(\Gamma_{3}),\Z),N\cap \Gamma_{3})\to Hom(H_{2}(H,\Z),N\cap \Gamma_{3})$ is trivial. Hence the desired central extension is trivial on $H$.\hfill $\Box$

As in Lemma~\ref{fpfpfp}, one proves that $\Gamma_{4}$ is of type ${\rm FP}_{\infty}$. Applying again the result of~\cite{bhms2002}, we get that $\Gamma_{4}$ is isomorphic to a direct product of surface groups, possibly with a free Abelian factor. To obtain the more precise statement of Theorem~\ref{fac-complet}, and to justify the fact $\Gamma_{4}$ does not contain any free non-Abelian factor, we argue as follows. We will call direct factor of the product
\begin{equation}\label{onemooore}
\Z^{k}\times \pi_{1}(\Sigma_{1})\times \cdots \times \pi_{1}(\Sigma_{r-k}).
\end{equation}
either one of the groups $\pi_{1}(\Sigma_{s})$ or one of the $\Z$ copy of $\Z^{k}$. The intersection of $\Gamma_{4}$ with each of the $r$ direct factors in~\eqref{onemooore} must be nontrivial because $\Gamma_{4}$ contains a copy of $\Z^{r}$. Indeed if one of these intersections was trivial, $\Gamma_{4}$ would embed into a direct product which does not contain $\Z^{r}$. We call $L_{1}, \ldots , L_{r}$ these intersections, where $L_{i} <\Z$ for $1\le i \le k$ and where $L_{i}<\pi_{1}(\Sigma_{i-k})$ for $i\ge k+1$. The proof of~\cite[p. 101]{bhms2002} shows that for $i\ge k+1$, $L_{i}$ is finitely generated and of finite index inside $\pi_{1}(\Sigma_{i-k})$. This implies that $\Gamma_{4}$ contains a finite index subgroup $\Gamma_{\ast}$ which is the direct product of each of its intersections with the factors in~\eqref{onemooore}.  The group $\Gamma_{\ast}$ now satisfies the conclusion of Theorem~\ref{fac-complet}.

\subsection{Cubulable K\"ahler manifolds}\label{sscm}

We now turn to the proof of Theorem~\ref{manifold}. Let $M$ be as in the statement of the theorem. Note that in particular, $M$ is aspherical. Applying Corollary~\ref{cckcc}, we see that a finite cover $M_{1}$ of $M$ has fundamental group isomorphic to a product of the form
\begin{equation}\label{eeeeeee}
\Z^{2l}\times \pi_{1}(S_{1})\times \cdots \times \pi_{1}(S_{m})
\end{equation}
where the $S_{i}$'s are closed surfaces of genus greater than $1$. We fix such an isomorphism and we denote by $\pi_{i}$ ($1\le i \le m$) the projection from $\pi_{1}(M_{1})$ onto $\pi_{1}(S_{i})$. From now on the proof will not make any further use of $\cat$ cubical complexes. We only use arguments from K\"ahler geometry. First, we will need the following classical result, see Theorem 5.14 in~\cite{catanese}:

\begin{theorem}  Let $X$ be a K\"ahler manifold, $S$ a topological surface of genus $\geqslant 2$, and $\pi : \pi_1(X)  \rightarrow \pi_1(S) $ a surjective homomorphism with finitely generated kernel. Then, there exists a complex structure on $S$ and a holomorphic map with connected fibers $p : X \rightarrow S$ such that the map
$$p_{\ast} : \pi_1(X)  \rightarrow \pi_1(S)$$
induced by $p$ is equal to $\pi$.
\end{theorem}

Applying this theorem to the various $\pi_{i}$'s we obtain that the surfaces $S_{i}$'s can be endowed with complex structures such that one has holomorphic maps $p_{i} : M_{1} \to S_{i}$ inducing the homomorphisms $\pi_{i}$ at the level of fundamental group. 

Let $A$ be the Albanese torus of $M_{1}$ and let $\alpha : M_{1} \to A$ be the Albanese map, which is well-defined up to translation. We also denote by $A_{i}$ the Albanese torus (or Jacobian) of the Riemann surface $S_{i}$ and by $\alpha_{i} : S_{i} \to A_{i}$ the corresponding map. By definition of the Albanese maps, for each $i$ there exists a holomorphic map $\varphi_{i} : A \to A_{i}$ which makes the following diagram commutative:
\begin{displaymath}
\xymatrix{M_{1} \ar[r]^{\alpha} \ar[d]^{p_{i}} & A \ar[d]^{\varphi_{i}} \\
S_{i} \ar[r]^{\alpha_{i}} & A_{i}. \\}
\end{displaymath} 
 We denote by $\varphi : A \to A_{1}\times \cdots \times A_{m}$ the product of the maps $\varphi_{i}$ and by $\beta : S_{1}\times \cdots \times S_{m}\to A_{1}\times \cdots \times A_{m}$ the product of the maps $\alpha_{i}$. Up to composing the maps $\alpha_{i}$'s with some translations, we can and do assume that $\varphi$ maps the origin of $A$ to the origin of $A_{1}\times \cdots \times A_{m}$, hence is a group homomorphism. 

Let $Y$ be the following submanifold of $A$:
$$Y=\{y\in A, \varphi (y)\in Im (\beta)\}$$
This is indeed a submanifold since $\beta$ is an embedding and $\varphi$ is a submersion. Now by construction the Albanese map $\alpha$ of $M_{1}$ can be written as:
$$\alpha = i\circ \Phi$$
where $i$ is the inclusion of $Y$ in $A$ and $\Phi: M_{1} \to Y$ is holomorphic. We now have:

\begin{lemma} The map $\Phi$ is a homotopy equivalence between $M_{1}$ and $Y$. 
\end{lemma} 
{\it Proof.} Let $B$ be the kernel of the map $\varphi$. The complex dimension of $B$ is equal to 
$$\frac{1}{2}\left(b_{1}(M_{1})-\sum_{j=1}^{m}b_{1}(S_{j})\right)$$
which equals the number $l$ appearing in Equation~\eqref{eeeeeee}. One can find a $C^{\infty}$ diffeomorphism $\theta : A\to B \times A_{1}\times \cdots \times A_{m}$ such that $\varphi \circ \theta^{-1}$ is equal to the natural projection 
$$ B \times A_{1}\times \cdots \times A_{m}\to  A_{1}\times \cdots \times A_{m}.$$
This implies that the complex manifold $Y$ is $C^{\infty}$ diffeomorphic to $B \times S_{1}\times \cdots \times S_{m}$. In particular, $Y$ and $M_{1}$ have isomorphic fundamental groups. We finally prove that $\Phi$ induces an isomorphism on fundamental groups. It follows from our description of $Y$ that one can choose identifications of $\pi_{1}(M_{1})$ and $\pi_{1}(Y)$ with 
$$\Z^{2l}\times \pi_{1}(S_{1})\times \cdots \times \pi_{1}(S_{m})$$
in such a way that $\Phi_{\ast}$ preserves each of the projections onto the groups $\pi_{1}(S_{j})$. This implies that $\Phi_{\ast}$ induces an isomorphism between the quotients of $\pi_{1}(M_{1})$ and $\pi_{1}(Y)$ by their respective centers. To conclude, it is enough to check that $\Phi_{\ast}$ induces an isomorphism between the centers of $\pi_{1}(M_{1})$ and $\pi_{1}(Y)$. But the composition of $\Phi_{\ast}$ with  the projection from $\pi_{1}(Y)$ onto its abelianization coincides with the map $\alpha_{\ast} : \pi_{1}(M_{1})\to \pi_{1}(A)\simeq H_{1}(M_{1},\Z)$. Since the center of $\pi_{1}(M_{1})$ injects in $H_{1}(M_{1},\Z)$, the restriction of $\Phi_{\ast}$ to the center of $\pi_{1}(M_{1})$ must be injective. Now using the fact that the quotient of $H_{1}(M_{1},\Z)$ by the image of the center of $\pi_{1}(M_{1})$ in $H_{1}(M_{1},\Z)$ is torsionfree, one sees easily that $\Phi_{\ast}(Z(\pi_{1}(M_{1})))$ must be equal to the center of $\pi_{1}(Y)$. This concludes the proof that $\Phi_{\ast}$ is an isomorphism. 

The manifold $M_{1}$ is aspherical by hypothesis. The manifold $Y$ is also aspherical since it is diffeomorphic to $B\times S_{1} \times \cdots \times S_{m}$. Since $\Phi$ induces an isomorphism on fundamental group, it is a homotopy equivalence. This concludes the proof of the lemma.\hfill $\Box$

We now conclude the proof using the following fact due to Siu and contained in the proof of Theorem~8 from~\cite{siu}:

\begin{center} {\it Let $f : Z_{1}\to Z_{2}$ be a holomorphic map between two compact K\"ahler manifolds of dimension $n$. Assume that $f$ is of degree $1$ and that the induced map $H_{2n-2}(Z_{1},\R)\to H_{2n-2}(Z_{2},\R)$ is injective. Then $f$ is a holomorphic diffeomorphism. The conclusion holds in particular if $f$ is a homotopy equivalence.} 
\end{center}

\noindent Applying Siu's result to the map $\Phi : M_{1} \to Y$ we obtain that $M_{1}$ and $Y$ are biholomorphic. This proves the first statement in Theorem~\ref{manifold}. When the original manifold $M$ is algebraic, an easy application of Poincar\'e's reducibility theorem~\cite[VI.8]{debarre} shows that a finite cover of $M_{1}$ is actually biholomorphic to a product of a torus with finitely many Riemann surfaces. This concludes the proof of Theorem~\ref{manifold}.

\section{Comments}\label{ques}

We discuss here some possible improvements to our results. 

First, one would like to remove the hypothesis of local finiteness in Theorem~\ref{fac-pv}. We summarize at which points this hypothesis was used:
\begin{enumerate}
\item It was used for the first time in Proposition~\ref{finitude}. However in this place, we have seen that it can be replaced by the hypothesis that the group action under consideration has finite stabilizers. Note that Proposition~\ref{finitude} is used later in section~\ref{sphm} to prove that a certain function defined as a supremum of continuous plurisubharmonic functions is actually the supremum of a finite number of continuous functions, hence is continuous. 
\item It is also used in the proof of Proposition~\ref{proprete} to show that the {\it signed distance function} to a hyperplane induces a proper function on a suitable covering space of the manifold under consideration. 
\end{enumerate} 
 
Second, one could try to remove the hypothesis that there is no fixed point in the visual boundary in Theorem~\ref{fac-pv}. There is a well-known strategy to achieve this, see~\cite[\S 2.H]{cif} as well as appendix B by Caprace in~\cite{cif}. If a group $G$ acts without fixed point on a $\cat$ cubical complex $X$ but with a fixed point in the visual boundary $\partial_{\infty}X$, one can construct another $\cat$ cubical complex $X_{1}$, of smaller dimension, on which $G$ acts and such that $X_{1}$ embeds in the Roller boundary of $X$. By applying this construction repeatedly, one can obtain a description of actions having a fixed point in the visual boundary. The reason why we cannot use this method here is that the passage from $X$ to $X_{1}$ does {\it not} preserve the local finiteness of the complex. We thank Pierre-Emmanuel Caprace for useful discussions concerning this point. Thus, one sees that removing the hypothesis of local finiteness from Theorem~\ref{fac-pv} should also allow to describe {\it parabolic} actions of K\"ahler groups on $\cat$ cubical complexes. Note that parabolic actions of K\"ahler groups on trees are already understood~\cite{delzant2}.


\bigskip
\bigskip
\begin{small}
\begin{tabular}{llll}
Thomas Delzant & & & Pierre Py\\
IRMA & & & Instituto de Matem\'aticas\\
Universit\'e de Strasbourg \& CNRS & & & Universidad Nacional Aut\'onoma de M\'exico\\
67084 Strasbourg, France & & & Ciudad Universitaria, 04510 M\'exico DF, M\'exico\\
delzant@math.unistra.fr & & & py@im.unam.mx\\    
\end{tabular}
\end{small}
 
 \end{document}